\documentclass[12pt]{article}

\usepackage[english]{babel}
\usepackage{amsmath,amsthm}
\usepackage{amsfonts}
\usepackage{graphicx}
\usepackage{epsfig}
\usepackage{epstopdf}
\usepackage{url}
\usepackage{amssymb}
\usepackage{color}
\usepackage{psfig}
\usepackage[%
breaklinks%
,colorlinks%
,linkcolor=red
,anchorcolor=blue%
,pagecolor=blue%
,citecolor=blue%
,bookmarks=false%
]{hyperref}

\theoremstyle{Definition}

\numberwithin{equation}{section}

\newtheorem{theorem}{Theorem}
\newtheorem{corollary}{Corollary}

\newtheorem{remark}{Remark}
\newtheorem{definition}{Definition}
\newtheorem{example}{Example}

\def\para{\vspace{4 mm}}

\def\cQ{{\mathcal Q}}
\def\cM{{\mathcal M}}
\def\para{\vspace{2.5mm}}
\def\cP{{\mathcal P}}

\def\gcd{{\rm gcd}}
\def\lcm{{\rm lcm}}

\def\deg{{\rm deg}}

\def\limit{{\lim}}

\begin{document}

 \title{Computing branches and asymptotes of meromorphic functions}

  \author{M. Fernández de Sevilla, R. Magdalena Benedicto$^*$ and S. P\'erez-D\'{\i}az\\
        Universidad de Alcal\'a and Universidad de Valencia$^*$ \\
 Dpto. de Ciencias de la Computación,\\  IDAL Electronic Engineering Department$^*$,   Dpto.  F\'{\i}sica y Matem\'aticas \\
      E-28871 Madrid and  E-46100 Burjassot Valencia$^*$ (Spain)  \\
 marian.fernandez@uah.es, rafael.magdalena@uv.es, sonia.perez@uah.es}

\date{}          
\maketitle

\begin{abstract}
In this paper, we first summarize the existing algorithms for computing all the {\it generalized
asymptotes} of a plane algebraic curve
implicitly or parametrically defined.  From these previous results, we derive  a  method that allows to easily compute the whole branch and all the {\it generalized asymptotes} of a ``special" curve   defined in  $n$-dimensional space by a parametrization that is not necessarily rational. So, some new concepts and methods are established for this type of curves. The approach is based on the notion of perfect curves introduced from
the concepts and results presented  in previous papers.
\end{abstract}


\noindent
{\it MSC Classification Codes:} 14H50, 30D30, 30E15,14A25

\section{Introduction}

An asymptote of a curve is a line to which the curve converges. In other words, the curve and its asymptote get infinitely close. Asymptotes have a variety of applications: they are used in big O notation, they are simple approximations to complex equations,  they are useful for graphing curves, etc.   Graphic means of displaying information are used in all areas of society. They have a complete image, are characterized by symbolism, compactness, relative ease of reading. It is these qualities of graphic images that determine their expanded use. In the near future, more than half of the information presented will have a graphical presentation form. The development of the theoretical foundations of descriptive geometry, engineering graphics, and other related sciences has expanded the methods for obtaining graphic images. Along with manual methods of forming graphic images, compiling project documentation, computer methods are finding wider application. The use of new information technologies provides the creation, editing, storage, replication of graphic images using various software tools.

In this sense the computation  of branches and asymptotes as a mathematical tool is very important since curves are essential  for engineering, industry, computer aided design (CAD), etc.

\para

For instance there are many applications of engineering curves in industry. The hyperbolic shape for example, finds application in design of cooling towers. Even Mirrors used in long telescopes are hyperbolic in shape. Another type of engineering curve called the archimedean spiral (type of curve in which the moving point is traced out in such a way that movement towards or away from the pole is uniform with vectorial angle from the starting line) has its application in designing and manufacturing of teeth profile of helical gears and profile of cams. Another curve is the cycloid which is used by engineers and designers for designing roller coasters. Even worm gears have cycloidal profile (the ones used in outdoor gears). The head of the tooth of such worm gear is an epicycloid (another engineering curve) and the tooth foot is hypocycloid. While designing objects various types of curves are used.

\para

Mechanical engineers also need mathematical curves. For example, a satellite dish is a basic parabola, a gear has the involute of a circle as its base. These kind of curves are usually not directly supported in CAD systems. They must therefore be drawn  using its branches.

\para

There are numerous methods for analysis and synthesis of mechanisms based on geometrical constructions and it is necessary a deepen study of the curves described by a point and the relationship between the geometry of different parts. Many engineering studies are devoted to the study of curves of the tooth profile of gears as well as the coupler path of mechanisms. Then, Geometry plays an important role in many engineering applications, such as engines and mechanisms. The study of curves dates from Ancient Greece, because the first mathematicians of History became interested in them. The Greeks were the first who studied the paths that describe planets in motion but they restricted their mathematics mainly to geometry, and they were primarily concerned with figures which could be obtained from lines and circles (geometric locus). Conics were treated as plane sections of cones (solid locus) and other planar curves like cycloids and spirals were included in their studies although they could not be drawn from lines and circles. Indeed they were known as mechanical curves rather than geometrical curves. In this paper we have focus the attention in drawing the mechanical curves most used in engineering by using dynamic geometry software; the different cycloid, hypocycloid, epicycloids have been drawn by using the Geogebra software. Some engineering applications of these mechanical curves, planetary gear trains, and the kinematic requirements have been also studied. For some bibliography see for instance \cite{Arn}, \cite{BaCr}, \cite{ChMar}, \cite{Kr}, \cite{Mar}, \cite{OletAl}, \cite{SoRe}, \cite{ZelMys}, \cite{Zwikker}, etc.

\para

The asymptotes of an infinity branch (a branch at infinity), $B$, of a real plane
algebraic curve, $\cal C$, reflect the behavior of $B$ at the points with
sufficiently large coordinates. In analytic geometry, an asymptote
of a curve is a line such that the distance between the curve and
the line approaches zero as they tend to infinity.  In some
contexts, such as algebraic geometry, an asymptote is defined as a
line which is tangent to a curve at infinity.

\para

If $B$ can be defined by some explicit equation of the form
$y = f(x)$ (or $x = g(y)$), where $f$ (or $g$) is a continuous function on an infinite interval, it is straightforward
to decide whether $\cal C$  has an asymptote at $B$ by analyzing the existence  of the limits
of certain functions when $x$ tends to $\infty$ (or $y$ tends to $\infty$). Moreover, if these limits  can be  computed, we may obtain the equation of the asymptote
of $\cal C$ at $B$. However, if this branch $B$ is implicitly defined and its equation cannot be converted into an
explicit form, both the decision and the computation of the asymptote of $\cal C$ at $B$ require some other tools. More precisely, an algebraic  curve may have more general curves than lines describing the behavior of a branch at the points
with sufficiently large coordinates.  Intuitively speaking, we say that a curve $\widetilde{{\cal C}}$ is a {\it  generalized asymptote} (or \textit{g-asymptote}) of another curve $\cal C$ if the distance between $\widetilde{{\cal C}}$ and $\cal C$  tends to zero as they tend to infinity, and $\cal C$  can not be approached by a new curve of lower degree (see  \cite{paper1}, \cite{paper2}, \cite{paper3} and \cite{newBP}). This motivates our interest in efficiently computing these {\it  generalized asymptotes} for a  wider variety of varieties such as the curves defined by a not necessarily rational parameterization.

\para

In this paper we deal with the problem of efficiently computing the asymptotes of meromorphic functions from an open subset of the complex plane onto ${\Bbb C}^n,\,n\geq 1$. We remind that  meromorphic functions are  functions on an open subset $D$ of $\Bbb C$ that are holomorphic on all  $D$ except for a set of isolated points, which are poles of the functions. By abuse of notation, and in order to  make the article easier for the reader to understand, we will denote by ${\cal P}(t):=(p_1(t),\ldots,p_n(t))$ these meromorphic functions and we say that the image of $\cal P$ is a {\it special curve $\cal C$ parametrically defined in   $n$-dimensional space.}

\para
The problem of the computation of asymptotes is dealt in previous papers of the third author (see \cite{paper1}, \cite{paper2}, \cite{paper3},  \cite{implem}, \cite{newBP}) and solved for algebraic rational curves parametrically and implicitly defined. For this purpose, some previous notions as infinity branches (or branches at infinity), approaching curves and perfect curves are introduced. The new goal we solve in this paper consists in working with curves parametrically defined but not necessarily rational.   This question is  very important in the study of these type of
curves because  there is no result or concept in this regard.

\para

We have intended the paper to be self-contained. For this reason, we have included Section \ref{S-notation}, where we review the theory of infinity branches and introduce the notions of convergent branches (that is, branches that get closer  as they tend to infinity) and approaching curves (see \cite{paper1}), and Section \ref{S-asymptotes}, where we lay down fundamental concepts like \textit{perfect curve} (a curve of degree $d$ that cannot
be approached by any curve of degree less than $d$) and \textit{g-asymptote} (a perfect curve that approaches another curve at an infinity branch).  In addition, we present the methods that allow to compute the infinity branches of a given curve implicitly and parametrically  defined, and  a g-asymptote for each of them (see Subsections \ref{Sub-asintotasconstructions}, \ref{Sub-asympparametriccurves} and \ref{Sub-improvement}).

\para

The main result of the paper is presented in Section \ref{S-new}. Here, we develop a  method that allows to easily compute all the generalized asymptotes of a   curve defined by a parametrization by only determining some simple limits of functions constructed from the given parametrization. The results presented are concerned with   plane curves but, as we remark in the paper, they can trivially be adapted for dealing with algebraic curves in   $n$-dimensional space (see Example \ref{ex-method-general-n-dimensional}).

\para
Finally,  some conclusions, and future work  is presented in  Section \ref{S-conclusion}.

\section{Notation and previous results}\label{S-notation}

In this section, we introduce the notion of {\it infinity branch or branch at infinity},  {\it convergent branches}
and  {\it approaching curves},
and we present some  properties which allow us to compare  the
behavior of two implicit algebraic plane curves at    infinity.
For more details on these concepts and results, we refer to
\cite{paper2} (see Sections 3 and 4).

\para

We consider an irreducible algebraic affine plane curve $\cal C$  over
$\mathbb{C}$  defined by the irreducible polynomial $f(x,y) \in {\Bbb R}[x,y]$.  We work over the field of complex numbers $\Bbb C$, but $\cal C$ has infinitely many points in the
affine plane over $\Bbb R$ (see Chapter 7 \cite{SWP}). The assumption of reality is included because
of the nature of the problem, but the theory can be similarly developed for the case of complex non-real curves.

Let ${\cal C}^*$  be its corresponding projective curve, defined by
the homogeneous polynomial $$F(x,y,z)=f_d(x,y)+zf_{d-1}(x,y)+z^2f_{d-2}(x,y)+\cdots+z^d f_0 \in {\Bbb R}[x,y,z],\quad $$
where $d:=\deg({\cal C})$ and $f_i(x,y),\,i=0,\ldots,d$ the homogeneous form of degree $i$. We assume that $(0:1:0)$ is not a  point at infinity of  ${\cal C}^*$ (otherwise, we may consider a projective linear change of coordinates).

\para

In order to get the infinity branches (or branches at infinity) of $\cal {C}$, we consider the  curve  defined by the
polynomial $g(y,z)=F(1:y:z)$ and we compute the
series expansion for the solutions of $g(y,z)=0$ around $z=0$.  We denote by  ${{\Bbb C}\ll z\gg}$ e the field of {\em  formal Puiseux series}. Thus, there exist exactly $\deg_{y}(g)$ solutions given by different Puiseux series
that can be grouped into  conjugacy classes.  More precisely, if $$\varphi(z)=m+a_1z^{N_1/N}+a_2z^{N_2/N}+a_3z^{N_3/N} +
\cdots\in{\Bbb C}\langle\langle z\rangle\rangle,\quad a_i\not=0,\, \forall i\in {\Bbb
N},$$
where $N\in {\Bbb
N}$ , $N_i\in {\Bbb
N},\,\,i\in {\Bbb N}$,  and $0<N_1<N_2<\cdots$,\,$\gcd(N,N_1,N_2,\ldots)=1$, is a Puiseux
series such that $g(\varphi(z), z)=0$, and   $\nu(\varphi)=N$  (i.e., $N$ is the  {\it  ramification index} of $\varphi$), the
 series
 \[\varphi_j(z)=m+a_1c_j^{N_1}z^{N_1/N}+a_2c_j^{N_2}z^{N_2/N}+a_3c_j^{N_3}z^{N_3/N}
+ \cdots\] where $c_j^N=1,\,\,j\in\{1,\ldots,N\}$,  are called
the {\it conjugates} of $\varphi$ (that is, $c_j,\,\,j\in\{1,\ldots,N\}$ are the   $N^{th}$  roots of unity).
The set  $\{\varphi_1,\ldots, \varphi_N\}$  all the conjugates of  $\varphi$ is called the  {\it  conjugacy
class} of $\varphi$. It contains $N=\nu(\varphi)$ distinct series which satisfy $g(\varphi_j(z), z)=0,\,j=1,\ldots,N$.

\para

Since
$g(\varphi(z), z)=0$ in some neighborhood of $z=0$ where
$\varphi(z)$ converges, there exists $M \in {\Bbb R}^+$
such that
$F(1:\varphi(t):t)=g(\varphi(t), t)=0$ for $t\in {\Bbb C}$ and
$|t|<M$, which implies that
$F(t^{-1}:t^{-1}\varphi(t):1)=f(t^{-1},t^{-1}\varphi(t))=0$, for
$t\in {\Bbb C}$  and  $0<|t|<M$. We set $t^{-1}=z$, and   we obtain
that $f(z,r(z))=0$ for $z\in {\Bbb C}$ and
$|z|>M^{-1}$ where
$$r(z)=z\varphi(z^{-1})=mz+a_1z^{1-N_1/N}+a_2z^{1-N_2/N}+a_3z^{1-N_3/N} + \cdots,\quad a_i\not=0,\, \forall i\in {\Bbb N} $$
$N,N_i\in {\Bbb N},\,\,i\in {\Bbb N}$, and $0<N_1<N_2<\cdots$,\,$\gcd(N,N_1,N_2,\ldots)=1$.

\para

Reasoning similarly with the $N$ different series in the conjugacy
class,  $\varphi_1,\ldots,\varphi_N$, we get
\[r_i(z)=z\varphi_i(z^{-1})=mz+a_1c_i^{N_1}z^{1-N_1/N}+a_2c_i^{N_2}z^{1-N_2/N}+a_3c_i^{N_3}z^{1-N_3/N}
+ \cdots.\]


\begin{definition}\label{D-infinitybranch}
 An {\it infinity branch (or branch at infinity) of an affine plane curve ${\cal C}$}  associated to the infinity
point (or point at infinity) $P=(1:m:0),\,m\in {\Bbb C}$, is  a set
$\displaystyle B=\bigcup_{j=1}^N L_j$, where $L_j=\{(z,r_j(z))\in
{\Bbb C}^2: \,z\in {\Bbb C},\,|z|>M\}$,\,  $M\in
{\Bbb R}^+$, and
\begin{equation}\label{E-branch}r_j(z)=z\varphi_j(z^{-1})=mz+a_1c_j^{N_1}z^{1-N_1/N}+a_2c_j^{N_2}z^{1-N_2/N}+a_3c_j^{N_3}z^{1-N_3/N}
+ \cdots\end{equation} where $N, N_i\in {\Bbb N},\,\,i\in {\Bbb N}$,  $0<N_1<N_2<\cdots$, and $c_j^N=1,\,\,j\in\{1,\ldots,N\}$.
 The subsets $L_1,\ldots,L_N$ are
called the {\it leaves} of the infinity branch $B$.
 \end{definition}

\para

\begin{remark} \label{R-conjugation} An infinity branch  is uniquely determined  from one leaf, up to conjugation.
\end{remark}

\para

By abuse of notation, in the following we write $B=\{(z,r(z))\in {\Bbb C}^2: \,z\in {\Bbb C},\,|z|>M\}$ (where $M:=\max\{M_1,\ldots,M_N\}$). We recall that $N$ is   the ramification index of the branch $B$ and we  will write
$N=\nu(B)$ (the branch $B$ has $\nu(B)$  leaves).

\para

\begin{remark} \label{R-infpoint} Each infinity branch is associated to
a unique infinity point. More precisely, as we stated above, there exists $M \in {\Bbb R}^+$
such that $F(1:\varphi(t):t)=g(\varphi(t), t)=0$ for $|t|<M$, where
$$\varphi(z)=m+a_1z^{N_1/N}+a_2z^{N_2/N}+a_3z^{N_3/N} +
\cdots\in{\Bbb C}\langle\langle z\rangle\rangle.$$
Thus, for $t=0$ we get the infinity point $P=(1:\varphi(0):0)=(1:m:0)\in {\cal C}^*.$

\para

Conversely, given an infinity point $P=(1:m:0)$, there must be, at least, one Puiseux solution $\varphi$ such that $\varphi(0)=m$; this solution provides an infinity branch associated to $P$.
In particular, we conclude that every algebraic plane curve has, at least, one infinity branch. 


\end{remark}

\para

The procedure introduced above allows us to obtain the infinity branches of a curve $\cal C$, under the assumption that $(0:1:0)\notin {\cal C}^*$. However, a curve may have infinity branches, associated to the infinity point $(0:1:0)$, which can not be constructed in this way. These infinity branches have the form $\{(r(z),z)\in {\Bbb C}^2: \,z\in {\Bbb C},\,|z|>M\}$ and may be obtained by interchanging the variables $x$ and $y$. See \cite{paper2} (Definition 3.3) for further details.

\para

In the following, we introduce the notions of convergent branches
and approaching curves. Intuitively speaking, two infinity branches
converge if they get closer  as they tend to infinity. This concept
will allow us to analyze whether two curves approach each other.

\para

\begin{definition}\label{D-distance0}
Two
infinity branches, $B$ and $\overline{B}$, are convergent if there
exist two leaves   $L=\{(z,r(z))\in {\Bbb C}^2:\,z\in {\Bbb C},\,|z|>M\} \subset B$ and
$\overline{L}=\{(z,\overline{r}(z))\in {\Bbb C}^2:\,z\in {\Bbb
C},\,|z|>\overline{M}\}\subset \overline{B}$  such that  $\lim_{z\rightarrow\infty} (\overline{r}(z)-r(z))=0.$ In this case, we say that the leaves $L$ and $\overline{L}$ converge.
\end{definition}

\para

The following theorem provides a characterization for the convergence of two infinity branches (see \cite{paper2}).

\para

\begin{theorem}\label{L-DistVertical} The following statements hold:
\begin{itemize}\item[1.] Two  leaves $L=\{(z,r(z))\in {\Bbb C}^2:\,z\in {\Bbb C},\,|z|>M\}$ and
$\overline{L}=\{(z,\overline{r}(z))\in {\Bbb C}^2:\,z\in {\Bbb
C},\,|z|>\overline{M}\}$ are convergent if and only if the terms
with non negative exponent in the series $r(z)$ and
$\overline{r}(z)$ are the same.
\item[2.]  Two infinity branches $B$ and $\overline{B}$ are convergent if and
only if for each leaf $L\subset B$ there exists a leaf $\overline{L}\subset
\overline{B}$ convergent with $L$, and conversely.
\item[3.]  Two convergent infinity branches must be associated to
the same infinity point.
\end{itemize}\end{theorem}

\para

This paper is concerned with the study of the asymptotes of a curve. The classical concept of asymptote stands for a line that \textit{approaches} a given curve when it tends to the infinity. In the following we generalize this idea by claiming that two curves \textit{approach each other} if they, respectively, have two infinity branches that converge.

\para

\begin{definition}\label{D-distance1}
Let ${\cal C} $ be an algebraic plane curve  with an
infinity branch $B$. We say that a  curve ${\overline{{\cal C}}}$
{\it approaches} ${\cal C}$ at its infinity branch $B$ if there
exists one leaf $L=\{(z,r(z))\in {\Bbb C}^2:\,z\in {\Bbb
C},\,|z|>M\}\subset B$ such that
$\lim_{z\rightarrow\infty}d((z,r(z)),\overline{\cal C})=0,$ where $d(\cdot,\cdot)$ represents the euclidean distance.
\end{definition}

\para
The following theorem characterize the convergence of two curves at an infinity branch (see \cite{paper2}).

\para

\begin{theorem}\label{T-curvas-aprox}
Let ${\cal C}$
be a plane algebraic curve with an infinity branch
$B$. A plane algebraic curve ${\overline{{\cal C}}}$ approaches
${\cal C}$ at $B$ if and only if ${\overline{{\cal C}}}$ has an
infinity branch, $\overline{B}$, such that $B$ and $\overline{B}$
are convergent.
\end{theorem}

\para

Obviously, ``approaching'' is a symmetric concept, that is, ${\cal C}_1$ approaches ${\cal C}_2$ if and only if ${\cal C}_2$ approaches ${\cal C}_1$. When it happens we say that ${\cal C}_1$ and ${\cal C}_2$ are \textit{approaching curves} or that they \textit{approach each other}. In the next section we use this concept to generalize the classical notion of asymptote of a curve.

\section{Asymptotes of an algebraic curve}\label{S-asymptotes}

Given an algebraic plane curve $\cal C$ and an infinity branch $B$, in Section \ref{S-notation}, we have described how $\cal C$ can  be
approached  at $B$ by a second curve ${\overline{\cal C}}$. Now, suppose that $\deg({\overline{\cal
C}})<\deg({\cal C})$. Then one may say that $\cal C$ {\it degenerates}, since it behaves at  infinity as a curve of
smaller degree. For instance, a hyperbola is a curve of degree 2 that has two real asymptotes, which implies that the hyperbola  degenerates, at
 infinity, to two lines. Similarly, one can check that every ellipse has two asymptotes, although they are complex lines in this case. However, the
asymptotic behavior of a parabola is different, since it cannot be approached at infinity by any line. This  motivates the
following definition:

\para

\begin{definition}\label{D-perfect-curve}
An algebraic curve of degree $d$ is a {\it perfect curve} if it cannot be approached by
any curve of degree less than $d$.
\end{definition}

\para

More properties on perfect curves can be found in \cite{paper1}. In particular, one has that if a given curve of degree $d$ has an only branch of degree $d$, then the input curve is perfect. For instance, a curve $\cal C$ defined by a proper parametrization of the form $(t^{n}, a_{n}t^{n}+a_{n-1}t^{n-1}+\cdots+a_0)$ is always perfect since it has an only branch $B$ given by $(z,r(z))=(z, a_{n}z+a_{n-1}z^{(n-1)/n}+\cdots+a_0)$ and $\deg({\cal C})=\deg(B)=n$ (see Definition \ref{D-degreebranch} for the degree of a branch).

\para

A curve that is not perfect can be approached by other curves of
smaller degree. If these curves are perfect, we call them
{\it g-asymptotes}. More precisely, we have the following definition.

\para

\begin{definition}\label{D-asymptote}
Let ${\cal C}$ be a curve with an infinity branch $B$. A {\it
g-asymptote} (generalized asymptote) of ${\cal C}$ at $B$ is a
perfect curve that approaches ${\cal C}$ at $B$.
\end{definition}

\para

The notion of {\it
g-asymptote} is similar to the classical concept of asymptote. The
difference is that a g-asymptote is not necessarily a line, but a
perfect curve. Actually, it is a generalization, since every line is
a perfect curve (this fact follows from Definition
\ref{D-perfect-curve}).  Throughout the paper we  refer sometimes to
{\it g-asymptote} simply as {\it asymptote}.

\para

\begin{remark}\label{R-minimal-degree}
The degree of a g-asymptote is less than or equal to the
degree of the curve it approaches. In fact, a g-asymptote of a curve
$\cal C$ at a branch $B$ has minimal degree among all the curves
that approach $\cal C$ at $B$.
\end{remark}

\para

 In Subsection \ref{Sub-asintotasconstructions}, we show that every infinity branch
of a given algebraic plane curve implicitly defined has, at least, one asymptote and we show how to compute it. For this purpose, we rewrite Equation \ref{E-branch}  defining a branch $B$ (see Definition \ref{D-infinitybranch}) as
 \begin{equation}\label{Eq-inf-branchn}r(z)=mz+a_1z^{1-n_1/n}+\cdots
+a_kz^{1-n_k/n}+a_{k+1}z^{1-N_{k+1}/N}+\cdots\end{equation} where $0<N_1<\cdots<N_k\leq N<N_{k+1}<\cdots$ and $\gcd(N,N_1,\ldots,N_k)=b$, $N=n \cdot b$, $N_j=n_j\cdot b,\,\,j\in\{1,\ldots,k\}$. That is, we have simplified the non negative
exponents such that $\gcd(n,n_1,\ldots,n_k)=1$. Note that $0<n_1<n_2<\cdots$, and  $n_k\leq n$, and $N<N_{k+1}$, i.e.  the terms $a_jz^{1-N_j/N}$
with $j\geq k+1$ are those which have negative exponent. We denote these terms as $A(z):=\sum_{\ell=k+1}^\infty
a_{\ell}z^{-q_{\ell}},$ where $q_{\ell}=1-N_{\ell}/N\in\mathbb{Q}^+,$ $\ell \geq k+1.$

\para

\noindent Under these conditions, we introduce the definition of degree of a branch $B$:

\para

\begin{definition}\label{D-degreebranch} Let $B=\{(z,r(z))\in
{\Bbb C}^2: \,z\in {\Bbb C},\,|z|>M\}$ ($r(z)$ is defined in (\ref{Eq-inf-branchn})) be an
infinity branch associated to an infinity point $P=(1:m:0), m\in\mathbb{C}$. We say that $n$ is the degree of $B$, and we denote it by $\deg(B)$.
\end{definition}

\subsection{Construction of a g-asymptote of a  curve implicitly defined}\label{Sub-asintotasconstructions}

Taking into account Theorems \ref{L-DistVertical} and \ref{T-curvas-aprox}, we have  that any
curve $\overline{{\cal C}}$ approaching ${\cal C}$ at $B$ should
have an infinity branch
 $\overline{B}=\{(z,\overline{r}(z))\in {\Bbb
C}^2:\,z\in {\Bbb C},\,|z|>\overline{M}\}$ such that the terms with
non negative exponent in $r(z)$ and $\overline{r}(z)$ are the same.\\

In the simplest case, if $A=0$ in the branch $B$ (i.e. there are no terms with negative
exponent; see equality (\ref{Eq-inf-branchn})), we could consider the branch
\begin{equation}\label{Eq-inf-branch3}
\tilde{r}(z)=mz+a_1z^{1-n_1/n}+a_2z^{1-n_2/n}+\cdots
+a_kz^{1-n_k/n},
\end{equation}
where $a_1,a_2,\ldots\in\mathbb{C}\setminus \{0\}$,\,$m\in {\Bbb
C}$, $n,n_1,n_2\ldots\in\mathbb{N}$, $\gcd(n,n_1,\ldots,n_k)=1$, and $0<n_1<n_2<\cdots$.
Note that $\tilde{r}$ has the same terms with non negative exponent as
$r$, and $\tilde{r}$ does not have terms with negative exponent.

\para

Let $\widetilde{{\cal C}}$ be the irreducible plane curve containing the branch
$\widetilde{B}=\{(z,\tilde{r}(z))\in {\Bbb C}^2:\,z\in {\Bbb
C},\,|z|>\widetilde{M}\}$ (note that $\widetilde{{\cal C}}$  is
unique since two different algebraic curves have finitely many
common points). Observe that
$$\widetilde{{\cal Q}}(t)=(t^n,mt^n+a_1t^{n-n_1}
+\cdots +a_kt^{n-n_k})\in {\Bbb C}[t]^2$$
is a polynomial  parametrization of $\widetilde{{\cal C}}$, and it is proper (see Lemma 3 in \cite{paper1}). In Theorem 2 in \cite{paper1}, we  prove that
$\widetilde{{\cal C}}$ is a  g-asymptote of ${\cal C}$ at $B$.

\para

From these results, we obtain the method presented in \cite{paper2} and \cite{paper3}, that computes   g-asymptotes and   that is independent of the leaf chosen to define  the infinity branch.  We assume that we have prepared the input curve $\cal C$, by means  of a suitable   projective linear change of coordinates, such that $(0:1:0)$ is not an infinity point of $\cal C$.

\para
\noindent
In the following, we illustrate the method with an example.

\para

\begin{example}\label{Ej-asint-impli} Let ${\cal C}$ be the curve of degree $d=6$  defined  by the irreducible polynomial\\

\noindent
$f(x,y)=8y^2x^3+16y^5x-42y^4x+164y^3x^2-71670yx-34853x^2+20428x+15075y^2x-2213y^3x+196x^3-2530y^2x^2+13946yx^2-56yx^3-4978y^3-15321y^2+175y^4-197y^5+8y^6 \in {\Bbb R}[x,y].$\\

\noindent
First, we have that $f_6(x,y)=y^5(y+2x)$. Hence,  the infinity points are
$P_1=(1:0:0)$ and $P_2=(1:-2:0).$\\

\noindent
We start by analyzing the point $P_1$: there  are three infinity branches associated to $P_1$, $B_{1j}=\{(z,r_{1j}(z))\in {\Bbb C}^2:\,z\in {\Bbb C},\,|z|>M_1\}$,\,$j=1,2,3$, where\\

\noindent
$ r_{11}(z)=7/2+7/2I+(-127/8-15/8I) z^{-1}+(2161/16+1189/32I) z^{-2}+  (-6553/4-84517/128I) z^{-3}+\cdots,$\\

\noindent
$r_{12}(z)=7/2-7/2I+(-127/8+15/8I) z^{-1}+(2161/16-1189/32I) z^{-2}+(-6553/4+84517/128I) z^{-3} +\cdots $\\

\noindent
$r_{13}(z)=5/24+ 10/3  z^{1/3}  {2}^{1/3}-   z^{2/3} {2}^{-1/3}+97/3 z^{-1}+ 1289/162  z^{-1/3} 2^{-4/3}+\cdots $\\

\noindent
(we compute $r_{1j},\,j=1,2,3$ using the {\sf algcurves} package included in  the computer algebra system {\sf Maple}; in particular we use the command {\sf puiseux}).\\
We compute $\tilde{r}_{1j}(z),\,j=1,2,3$, and we have that
$$\tilde{r}_{11}(z)=7/2+7/2I,\,\quad \tilde{r}_{12}(z)=7/2-7/2I,\,\quad \tilde{r}_{13}(z)=5/24+10/3 \cdot 2^{1/3}z^{1/3}-2^{-1/3}z^{2/3}.$$
The parametrizations of the asymptotes $\widetilde{\cal C}_j,\,j=1,2,3,$ are given by
$$\widetilde{\cal Q}_1(t)=(t,\,7/2+7/2I),\,\quad \widetilde{\cal Q}_2(t)=(t,\,  7/2-7/2I),\,\quad \widetilde{\cal Q}_3(t)=(t^3,\,  5/24+10/3 \cdot 2^{1/3}t-2^{-1/3}t^2)$$
which define two complex lines and the curve defined by the implicit polynomial  $$\tilde{f}_3(x,y)=-6912x^2-138240yx+1052800x+8640y^2-1800y+125-13824y^3\in {\Bbb R}[x,y] $$
(one may compute the polynomial defining implicitly $\widetilde{\cal C}_3$ using  for instance the results in \cite{SWP}; see Chapter 4).\\

Now, we focus on the point $P_2$: there  one infinity branch associated to $P_2$,  $B_{2}=\{(z,r_{2}(z))\in {\Bbb C}^2:\,z\in {\Bbb C},\,|z|>M_2\}$, where\\

\noindent
$r_{2}(z)=17-2z-261/4 z^{-1}-2241/8 z^{-2}+\cdots.$\\

\noindent
We compute $\tilde{r}_{2}(z)$, and we have that
$$\tilde{r}_{2}(z)=-2z+17.$$
The parametrization  of the asymptote  $\widetilde{\cal C}_4$ is given by
$$\widetilde{\cal Q}_4(t)=(t,\,-2 t+17) $$
that defines a line implicitly defined by the polynomial
$$\tilde{f}_4(x,y)=-2 x+17-y\in {\Bbb R}[x,y].$$

 In Figure \ref{F-ex1}, we plot the curve
$\cal C$, and the asymptotes $\widetilde{\cal C}_3$ and $\widetilde{\cal C}_4$ (the  asymptotes $\widetilde{\cal C}_1$ and $\widetilde{\cal C}_2$ are complex lines).

\begin{figure}[h]
$$
\begin{array}{cc}
\psfig{figure=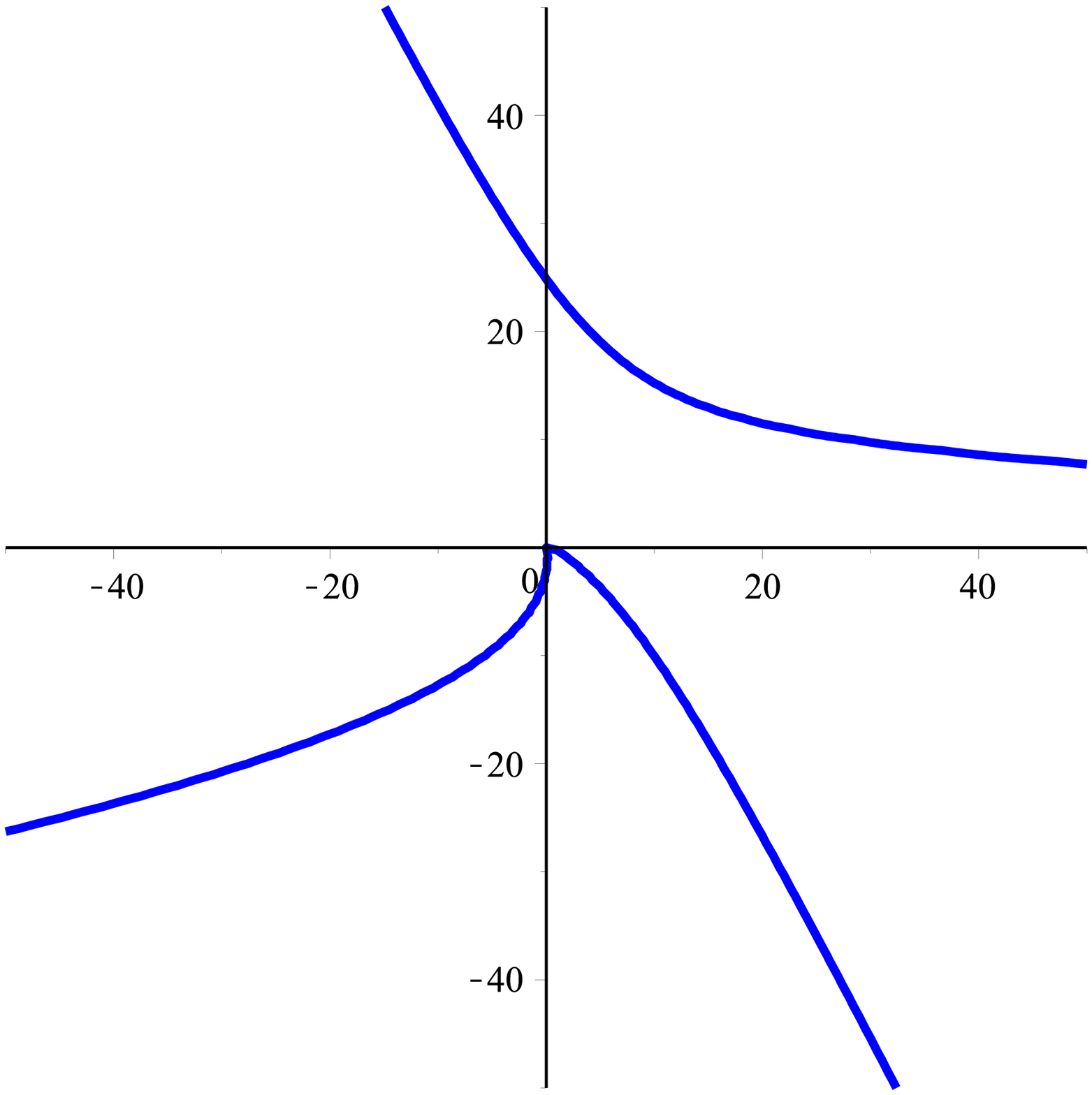,width=6 cm,height=6 cm} &
\psfig{figure=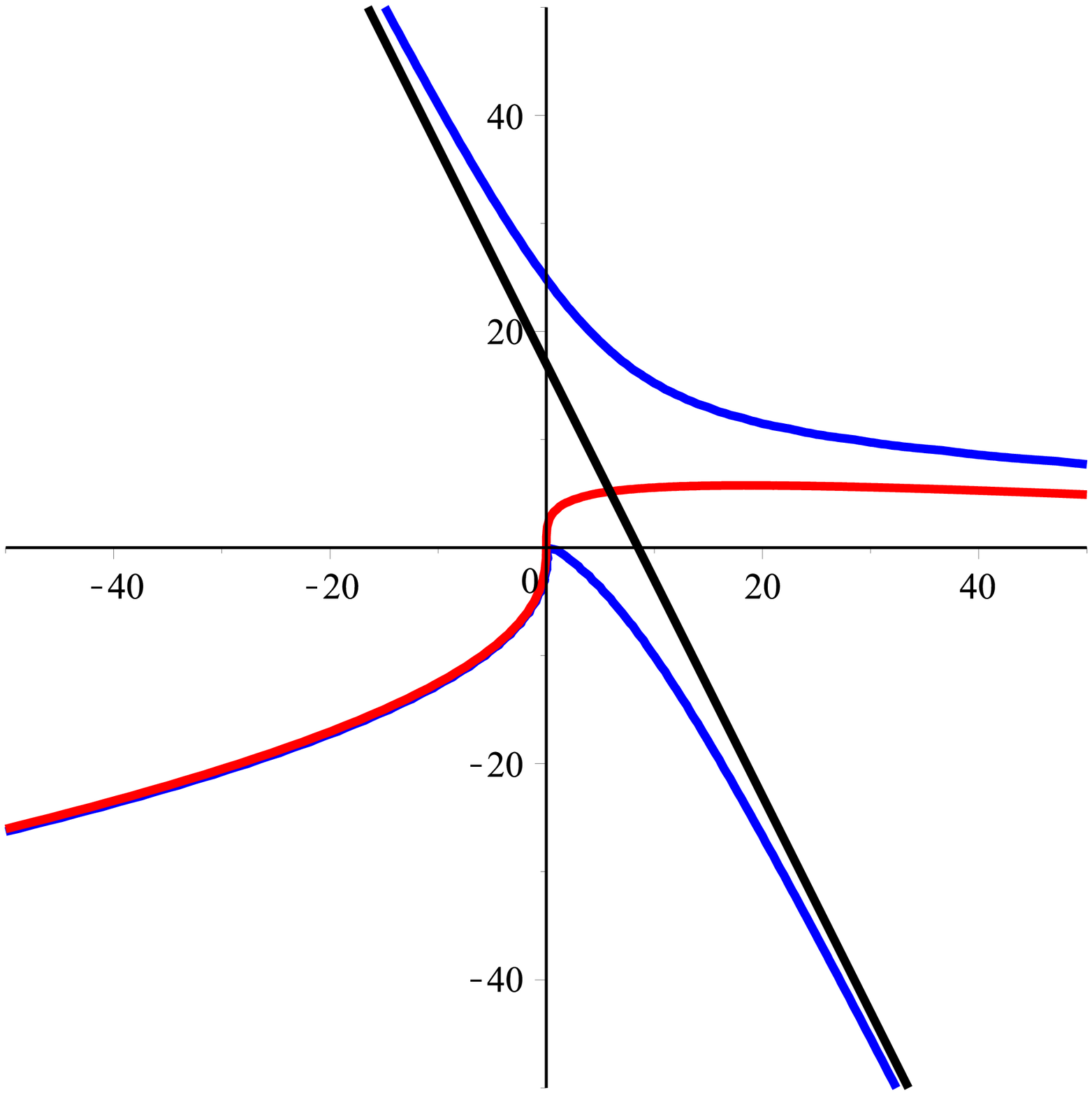,width=6 cm,height=6 cm}
\end{array}
$$ \caption{Curve $\cal C$ (left) and curve and asymptotes (right).}\label{F-ex1}
\end{figure}
\end{example}

\subsection{Construction of a g-asymptote of a curve rationally parametrized}\label{Sub-asympparametriccurves}

Throughout this paper so far, we have dealt with  algebraic plane curves
implicitly defined. In this subsection, we present a method
to compute infinity branches and g-asymptotes of a plane curve  from
their parametric (rational) representation (without implicitizing). This method is included in \cite{paper3} (see Section 5) and it involves the computation of Puiseux series and infinity branches. In Subsection \ref{Sub-improvement}, we  develop a new  method presented in \cite{newBP} that allows to easily compute the generalized asymptotes ({\it g-asymptotes}) by only determining some simple limits of rational functions constructed from the given parametrization.

\para

\noindent
Let ${\cal C}$  be a plane curve defined by the rational parametrization
\[{\cal P}(s)=(p_1(s),p_2(s))\in {\Bbb R}(s)^2,\quad p_i(s)=p_{i1}(s)/p_{i2}(s),\quad \gcd(p_{i1},p_{i2})=1,\,\,i=1,2.\]
If ${\cal C}^*$ represents  the projective curve associated to
${\cal C}$, we have that a parametrization of ${\cal C}^*$ is given by
${\cal P}^*(s)=(p_{1}(s):p_{2}(s):1)$ or, equivalently,
$${\cal P}^*(s)=\left(1:\frac{p_{2}(s)}{p_{1}(s)}:\frac{1}{p_{1}(s)}\right).$$
We assume that   we have  prepared the input curve ${\cal C}$, by means
of a suitable    projective linear change of coordinates (if necessary) such that
$(0:1:0)$ is not a point  at infinity of
${\cal C}^*$.

\para

In order to compute the g-asymptotes of $\cal C$, first we need to
determine the  infinity branches of $\cal C$. That is, the sets
$$B=\{(z, r(z))\in {\Bbb C}^2: \,z\in {\Bbb
C},\,|z|>M\},\,\, \mbox{where}\,\,
r(z)=z\varphi(z^{-1}).$$

 For this purpose, taking into account Definition \ref{D-infinitybranch}, we have that
$f(z,r(z))=F(1:\varphi(z^{-1}):z^{-1})=F(1:\varphi(t):t)=0$
around $t=0$,  where  $t=z^{-1}$ and  $F$ is the
polynomial  defining implicitly ${\cal C}^*$.  Observe that in this
section, we are given the parametrization ${\cal P}^*$ of ${\cal
C}^*$ and  then,  $F({\cal
P}^*(s))=F\left(1:{p_{2}(s)}/{p_{1}(s)}:{1}/{p_{1}(s)}\right)=0.$
Thus, intuitively speaking, in order to compute the  infinity
branches of $\cal C$, and in particular the series
$\varphi$, one needs to rewrite the parametrization
${\cal
P}^*(s)$
in the form $(1:\varphi(t):t)$ around $t=0$. For
this purpose, the idea is to look for a value of the parameter $s$,
say $\ell(t)\in {\Bbb C}\langle\langle
t\rangle\rangle $, such that ${\cal P}^*(\ell(t))=(1:\varphi(t):t)$ around $t=0$.\\

Hence, from the above reasoning, we deduce that first, we have to consider the equation $1/p_{1}(s)=t$ (or equivalently,  $p_{12}(s)-tp_{11}(s)=0$), and we solve it in the variable $s$ around $t=0$. From Puiseux's Theorem,   there
exist solutions $\ell_1(t),\ell_2(t),\ldots,\ell_k(t)\in {\Bbb C}\langle\langle
t\rangle\rangle $ such that,
$p_{12}(\ell_i(t))-tp_{11}(\ell_i(t))=0,\,i\in\{1,\ldots,k\},$ in a neighborhood of $t=0$.\\

Thus,   for each $i\in\{1,\ldots,k\}$, there exists  $M_i\in {\Bbb R}^+$ such that
the points $(1:\varphi_{i}(t):t)$ or  equivalently, the
points $(t^{-1}:t^{-1}\varphi_{i}(t):1)$, where
$\varphi_{i}(t)=\frac{p_{2}(\ell_i(t))}{p_{1}(\ell_i(t))},$  are in ${\cal C}^*$ for $|t|<M_i$ (note that ${\cal P}^*(\ell(t))\in {\cal C}^*$ since ${\cal P}^*$ is a parametrization of ${\cal C}^*$). Observe  that $\varphi_{i}(t)$ is
a Puiseux series, since $p_{2}(\ell_i(t))$ and
$p_{1}(\ell_i(t))$ can be written as Puiseux series and ${\Bbb C}\langle\langle
t\rangle\rangle $ is a
field.\\

Finally, we set $z=t^{-1}$. Then, we have that the points
$(z,r_{i}(z))$, where $r_{i}(z)=z\varphi_{i}(z^{-1})$, are in ${\cal C}$ for $|z|>M_i^{-1}$. Hence, the infinity
branches of $\cal C$ are the sets
$B_i=\{(z,r_{i}(z))\in {\Bbb C}^2: \,z\in {\Bbb
C},\,|z|>M_i^{-1}\},\quad i\in\{1,\ldots,k\}.$

\para

Note that the series $\ell_i(t)$ satisfies that
$p_1(\ell_i(t))t=1$, for $i\in\{1,\ldots,k\}$. Then, we have that
$$\varphi_{i}(t)=\frac{p_{2}(\ell_i(t))}{p_{1}(\ell_i(t))}=p_2(\ell_i(t))t,\quad r_{i}(z)=z\varphi_{i}(z^{-1})=p_2(\ell_i(z^{-1})).$$


\para
Once we have the infinity branches, we can compute a g-asymptote for
each of them by simply removing the terms with negative exponent
from $r_{i}$.
\para

\para

Additionally we note, that   some of the  solutions $\ell_1(t),\ell_2(t),\ldots,\ell_k(t)\in
{\Bbb C}\langle\langle
t\rangle\rangle $ might belong to the same conjugacy class. Thus, we only consider one
solution for each of these classes. The output asymptote $\widetilde{\cal C}$   is independent of the solutions $\ell_1(t),\ell_2(t),\ldots,\ell_k(t)\in
{\Bbb C}\langle\langle
t\rangle\rangle $ chosen in step $1$, and of the leaf chosen to define  the branch $B$.

\para

In the following example, we consider a parametric plane curve with two
real infinity branches. We  obtain these branches and compute a g-asymptote for each of them.

\para

\begin{example}\label{E-method-ant}
The plane curve  ${\cal C}$  introduced in Example \ref{Ej-asint-impli}  turns out to be rational, parametrized by
$${\cal P}(s)=\left(\frac{3s^4-s-4+5s^3}{(s-1)s^3(s^2+1)}, \frac{2s^2-7s+2}{(s-1)s^2}\right)\in {\Bbb R}(s)^2.$$
We  compute the asymptotes of $\cal C$. For this purpose, we determine the solutions of the equation
$p_{12}(s)-tp_{11}(s)=0$
around $t=0$. For this purpose, we may use, for
instance, the command {\sf puiseux} included in   the  package {\sf
algcurves} of  the computer algebra system {\sf Maple}. There are
four solutions (up to conjugation) that are given by the Puiseux series
$$\ell_1(t)=1+411/16 t^3+21/2 t^2+3/2+\cdots,$$$$\ell_2(t)=-315/32 t^3+11279/64 I t^3+13/4 t^2-243/16 I t^2-5/4 t+7/4 I t+I+\cdots $$
$$\ell_3(t)=- 315/32 t^3- 11279/64 I t^3+ 13/4 t^2+ 243/16 I t^2- 5/4 t- 7/4 I t-I+\cdots,$$$$ \ell_4(t)=- 17/3 t^2- 42179/7776 2^{1/3}t^{5/3}- 2105/1296 2^{2/3}t^{4/3}+ 1/3 t+ 5/6 2^{1/3}t^{2/3}+2^{2/3}t^{1/3}+\cdots. $$

\noindent
Now, we compute\\

\noindent $r_{1}(z)=p_2(\ell_1(z^{-1}))= -2 z+17-261/4 z^{-1}+\cdots$\\

\noindent $r_{2}(z)=p_2(\ell_2(z^{-1}))=7/2-7/2 I-127/8 z^{-1} + 15/8 I z^{-1}+\cdots $\\

\noindent $r_{3}(z)=p_2(\ell_2(z^{-1}))=7/2+7/2 I-127/8 z^{-1} - 15/8 I z^{-1}+\cdots $\\

\noindent $r_{4}(z)=p_2(\ell_2(z^{-1}))=5/24-  2^{-1/3} z^{2/3}+ 10/3\cdot  2^{1/3}z^{1/3}+ 1289/648\cdot 2^{2/3} z^{-1/3}+ 173813/15552\cdot 2^{1/3} z^{-2/3}+\cdots $\\

\noindent
(we may use, for instance, the command {\sf series}
included in   the computer algebra system {\sf Maple}). The
curve has four infinity branches given by
$B_i=\{(z,r_{i}(z))\in {\Bbb C}^2: \,z\in {\Bbb
C},\,|z|>M\}$
for some $M\in {\Bbb R}^+$ (note that $B_4$   has three leaves).

\para

We obtain $\tilde{r}_{i}(z)$ by
removing the terms with negative exponent in $r_{i}(z)$ for $i=1,2,3,4$. We get
$$\tilde{r}_{1}(z)=-2 z+17\quad\text{ and }\quad\tilde{r}_{2}(z)= 7/2-7/2 I$$
$$\tilde{r}_{3}(z)= 7/2+7/2 I\quad\text{ and }\quad\tilde{r}_{4}(z)= 5/24-  2^{-1/3} z^{2/3}+ 10/3\cdot  2^{1/3}z^{1/3}.$$

The input curve $\cal C$ has two complex asymptotes $\widetilde{{\cal C}}_i$ at
$B_i$ for $i=1,2$ and two real  asymptotes $\widetilde{{\cal C}}_i$ at
$B_i$ for $i=3,4$  that can be polynomially parametrized by (see Figure \ref{F-ex1}):
$$\widetilde{\cal Q}_1(t)=(t,\,  -2t+17),\quad \widetilde{\cal Q}_2(t)=(t,\,  7/2-7/2 I),\quad \widetilde{\cal Q}_3(t)=(t,\,  7/2+7/2 I) $$
 $$\widetilde{\cal Q}_4(t)=(t^3,\,  5/24+10/3 \cdot 2^{1/3}t-2^{-1/3}t^2). $$
Compare the output with the output obtained in Example \ref{Ej-asint-impli}.
\end{example}

\para

\begin{remark}\label{R-param-l}\ \begin{enumerate}
                 \item When we compute the series ${\ell}_i$, we cannot handle its infinite terms so it must be truncated,
which may distort the computation of the series $r_{i}$. However, this distortion may not affect to all the terms in
$r_{i}$. In fact, the number of affected terms depends on the
number of terms considered in ${\ell}_i$. Nevertheless, note that we do not
need to know the full expression of $r_{i}$ but only the terms
with non negative exponent. In \cite{paper3} (Proposition 2), it is proved that one can get the
terms with non negative exponent in  $r_{i}$ by considering just $2\deg(p_1)+1$ terms of $\ell_i$.
                 \item We remind that before to apply the method,  the input curve must be prepared such that $(0:1:0)$ is not a point  at infinity. As an alternative, one could apply the algorithm first for $\cP$ and then for ${\overline{\mathcal{P}}}:=(p_2(s),p_1(s))\in {\Bbb R}(s)^2$. In this last case, if we get the asymptote $(h_1,\,h_2)$, we have to undo the necessary change of coordinates and we finally get the asymptote $\widetilde{\cal Q}(t)=(h_2,\,h_1)$. Some of the asymptotes obtained from $\overline{\mathcal{P}}$ may coincide with others obtained from $\cP$ but some other new asymptotes could appear (those corresponding to vertical asymptotes; see Corollary \ref{C-final}.
               \end{enumerate}

\end{remark}

\subsubsection{New  method for the parametric (rational) case}\label{Sub-improvement}

In this subsection, we present an improvement of the method described above, which  avoids the computation of infinity branches and Puiseux series (see \cite{newBP}). We develop this method for the plane case but it can be trivially adapted  for dealing with rational curves in   $n$-dimensional space.

\para

\noindent In the following we consider a rational plane curve $\cal C$ defined  by the rational parametrization
\[{\cal P}(s)=(p_1(s),p_2(s))\in {\Bbb R}(s)^2,\quad p_i(s)=p_{i1}(s)/p_{i2}(s),\quad \gcd(p_{i1},p_{i2})=1,\,\,i=1,2.\]

We assume that $\deg(p_{i1})\leq \deg(p_{i2})=d_i,\,i=1,2$ (otherwise, we apply a suitable linear change on the variable  $t$). Thus, we have that
$\lim_{s\rightarrow\infty}p_i(s)\neq\infty,\,i=1,2$ and the infinity branches of $\cal C$ will be traced when $s$ moves around the different roots of the denominators $p_{12}(s)$ and $p_{22}(s)$.
In fact, each of these roots yields an infinity branch. The following theorem shows how to obtain a g-asymptote for each of these branches, by just computing some simple limits of rational functions constructed from  $\cP(s)$ (see \cite{newBP}).
\para

\begin{theorem}\label{T-final} Let  $\cal C$ be a curve defined  by a parametrization
\[{\cal P}(s)=(p_1(s),p_2(s))\in {\Bbb R}(s)^2,\quad p_i(s)=p_{i1}(s)/p_{i2}(s),\quad \gcd(p_{i1},p_{i2})=1,\,\,i=1,2,\]
where $\deg(p_{i1})\leq \deg(p_{i2})=d_i,\,i=1,2$. Let $\tau\in {\Bbb C}$ be such that $p_{i2}(t)=(t-\tau)^{n_i}\overline{p}_{i2}(t)$  (that is, $\tau\in {\Bbb C}$ is a root of multiplicity $n_i$ of $p_{i2}$) where $\overline{p}_{i2}(\tau)\not=0,\,i=1,2$ (that is, $\tau\in {\Bbb C}$ is not a root of $\overline{p}_{i2}$), and $n_1\geq 1$, and let $B$ be the corresponding infinity branch. A g-asymptote of $B$ is defined by the parametrization
$$\widetilde{\cal Q}(t)=(t^{n_1},\,  a_{n_2}t^{n_2}+a_{n_2-1}t^{n_2-1}+\ldots+a_{0}),$$
where
$$\begin{array}{ll}
   a_{n_2}=\limit_{t\rightarrow \tau}\,\, \displaystyle\frac{p_2(t)}{p_{1}(t)^{n_2/n_1}}& \\
  a_{n_2-1}=\limit_{t\rightarrow \tau}\, p_1(t)^{1/n_1}f_1(t),\quad &f_1(t):=\displaystyle \frac{p_2(t)}{p_{1}(t)^{n_2/n_1}}-a_{n_2}  \\
 a_{n_2-2}=\limit_{t\rightarrow \tau}\, p_1(t)^{1/n_1}f_2(t),\quad &f_2(t):= p_1(t)^{1/n_1}f_1(t)-a_{n_2-1}  \\
\qquad\quad \vdots & \qquad\quad \vdots\\
   a_{n_2-i}=\limit_{t\rightarrow \tau}\, p_1(t)^{1/n_1}f_i(t),\quad &f_i(t):= p_1(t)^{1/n_1}f_{i-1}(t)-a_{n_2-(i-1)},\,\,i\in\{2,\ldots,n_2\}. \\
 \end{array}$$

\end{theorem}

\para

\begin{remark}\label{R-infpoint2}
From the above construction, each root $\tau$ of $p_{12}(t)$ yields an infinity branch and, hence, an infinity point $P^*$ (see Remark \ref{R-infpoint}). Note that the parametrization $\mathcal{P}(t)$ can be expressed as
$\mathcal{P}(t)=\left(\frac{q_{11}(t)}{q(t)},\frac{q_{12}(t)}{q(t)}\right)$,
where $q(t)=\lcm(p_{12}(t),p_{22}(t))$ and $q_{1i}(t)=p_i(t)q(t)$. Now, the corresponding projective curve is parametrized by $\mathcal{P}^*(t)=(q_{11}(t),q_{12}(t),q(t))$
and the infinity point associated to $\tau$ is $P^*=(q_{11}(\tau):q_{12}(\tau):0)$.

\end{remark}

\para

In the following corollary, we analyze the special case of the vertical and horizontal g-asymptotes, i.e. lines of the form $x-a$ or $y-b$, where $a,b\in {\Bbb C}$ (observe that these asymptotes correspond to branches associated to the infinity points $(0:1:0)$ and $(1:0:0)$, respectively). More precisely, we prove that these asymptotes are obtained from the non--common roots of the denominators of the given parametrization. Note that  in the practical design  of engineering and modeling applications, the rational curves  are usually presented by numerical coefficients and $\cP(s)$    mostly satisfies that $\gcd(p_{12}, p_{22})=1$.

\para

\begin{corollary}\label{C-final}  Let  $\cal C$ be a curve defined  by a parametrization
\[{\cal P}(s)=(p_1(s),p_2(s))\in {\Bbb R}(s)^2,\quad p_i(s)=p_{i1}(s)/p_{i2}(s),\quad \gcd(p_{i1},p_{i2})=1,\,\,i=1,2,\]
where $\deg(p_{i1})\leq \deg(p_{i2}),\,i=1,2$.
\begin{enumerate}
\item Let $\tau\in {\Bbb C}$ be such that $p_{12}(t)=(t-\tau)^{n_1}\overline{p}_{12}(t)$ where $p_{22}(\tau)\overline{p}_{12}(\tau)\not=0$, and $n_1\geq 1$. It holds that a g-asymptote of $\cal C$ corresponding to the infinity point $(1:0:0)$ is the horizontal line  $y-p_{2}(\tau)=0$, defined by the parametrization
$\widetilde{\cal Q}(t)=(t,\,  p_{2}(\tau)).$
\item  Let $\tau\in {\Bbb C}$ be such that $p_{22}(t)=(t-\tau)^{n_2}\overline{p}_{22}(t)$ where $p_{12}(\tau)\overline{p}_{22}(\tau)\not=0$, and $n_2\geq 1$. It holds that a g-asymptote of $\cal C$ corresponding to the infinity point $(0:1:0)$ is the vertical line  $x-p_{1}(\tau)=0$, defined by the parametrization
$\widetilde{\cal Q}(t)=(p_{1}(\tau),\,  t).$
\end{enumerate}
   \end{corollary}

\para

\begin{remark}\label{R-proper} The previous theorem outputs the parametrization
 $\widetilde{\cal Q}(t)=(t^{n_1},\,  a_{n_2}t^{n_2}+a_{n_2-1}t^{n_2-1}+\ldots+a_{0}),$
and $n_1\geq n_2$ (otherwise $(0:1:0)$ is an infinity point of the input curve). Note that the degree of the defined curve is not necessary $n_1$ since $\cQ$ could be improper which is equivalent to $\gcd(n_1, n_2,\ldots,  n_2-j)\not=0$ for every $j=0,\ldots, n_2-1$ such that $a_{n_2-j}\not=0$.
  Let us assume that $\gcd(n_1, n_2,\ldots,  n_2-j)=\beta$ for every  $j=0,\ldots, n_2-1$ such that $a_{n_2-j}\not=0$. Then, let $n=n_1/\beta$ and
  \[\cM(t)=\cP(t^{1/\beta})=(t^{n}, a_{n_2}t^{n_2/\beta}+a_{n_2-1}t^{(n_2-1)/\beta}+\ldots+a_{0})\in {\Bbb K}[t]^2\]
is a proper reparametrization of $\cQ$. Then we get that the theorem outputs an asymptote since the output curve is perfect (it has an only branch and the degree of the curve which is  $n$ is equal to the degree of the branch).
   \end{remark}

\para

By applying the above results, we can easily obtain all the g-asymptotes of any rational plane curve, as  the following example  shows.

\para

\begin{example}\label{E-method-ant-b} Let ${\cal C}$ be the plane curve introduced in Examples \ref{Ej-asint-impli} and \ref{E-method-ant} defined by the
parametrization
$${\cal P}(s)=\left(\frac{3s^4-s-4+5s^3}{(s-1)s^3(s^2+1)}, \frac{2s^2-7s+2}{(s-1)s^2}\right)\in {\Bbb R}(s)^2.$$
We  compute the asymptotes of $\cal C$ using the new method just presented. For this purpose, we first  observe that $p_{12}(s)$ has the roots $\tau_1=1,\,\tau_2=0,\,\tau_3=I,\,\tau_4=-I$, with multiplicities $n_{11}=1$,   $n_{12}=3$ and $n_{13}=n_{14}=1$. The multiplicities of these roots in $p_{22}(s)$ are $n_{21}=1$,  $n_{22}=2$  and $n_{23}=n_{24}=0$.

\para

\noindent
 For $\tau_1=1,$ we compute
$$   a_{1}=\limit_{t\rightarrow 1}\,\frac{p_2(t)}{p_{1}(t)}=-2,\qquad
  a_{0}=\limit_{t\rightarrow 1}\, p_1(t)f_1(t)=17,\quad  f_1(t):=\frac{p_2(t)}{p_{1}(t)}-a_{1}.$$
Then, we obtain the asymptote $\widetilde{\cal C}_1$, defined by the proper parametrization
$$\widetilde{\cal Q}_1(t)=(t,\,  -2t+17).$$

\noindent
For $\tau_2=0,$ we compute
$$\begin{array}{ll}
   a_{2}=\limit_{t\rightarrow 2}\,\frac{p_2(t)}{p_{1}(t)^{2}}=2^{-4/3}+I\cdot 2^{-4/3}\,3^{1/2} &\\
  a_{1}=\limit_{t\rightarrow 2}\, p_1(t)f_1(t)=- 5/3\cdot 2^{1/3}+5/3\cdot I\cdot 2^{1/3}\,3^{1/2},\quad &f_1(t):=\frac{p_2(t)}{p_{1}(t)^{2}}-a_{2} \\
 a_{0}=\limit_{t\rightarrow 2}\, p_1(t) f_2(t)=\frac{5}{24},\quad &f_2(t):=p_1(t)f_{1}(t)-a_{1}. \\
 \end{array}$$
Then, we obtain the asymptote $\widetilde{\cal C}_2$, defined by the proper parametrization
$$\widetilde{\cal Q}_2(t)=\left(t^3,\,  5/24+(-5/3\cdot 2^{1/3}+5/3\cdot I\cdot 2^{1/3}\,3^{1/2})t+(2^{-4/3}+I\cdot 2^{-4/3}\,3^{1/2})t^2\right).$$

Finally for $\tau_3=I,$ and $\tau=-I$ we get the asymptotes $\widetilde{\cal C}_i,\,i=3,4,$, defined by the proper parametrizations
$$\widetilde{\cal Q}_3(t)=\left(t,\,  p_2(I)\right)=\left(t,\,  7/2-7/2 I\right),\quad  \widetilde{\cal Q}_4(t)=\left(t,\,  p_2(-I)\right)=\left(t,\,  7/2+7/2 I\right).$$
See Figure \ref{F-ex1} and  compare the output with the output obtained in Example \ref{E-method-ant}. We may check that $\widetilde{\cal Q}_2(t)$ is a reparametrización of the proper parametrization $\widetilde{\cal Q}_4(t)$ obtained in  Example \ref{E-method-ant}. In fact, note that  $\widetilde{\cal Q}_2(t)= \widetilde{\cal Q}_4(\xi\,t)$, with $\xi=1/2+I\sqrt{3}/2$ satisfies that $\xi^3=1$.
\end{example}

\para

\begin{remark}
 The above method allows us to easily obtain all the generalized asymptotes of a rational curve. However,  we should compute the roots of the denominators of the parametrization, which may entail certain difficulties if algebraic numbers are involved. This problem is solved using the notion of {\it conjugate points} (see Definition 12 in \cite{Perez}), which   help us to overcome this problem. The idea is to collect the points whose coordinates depend algebraically on all the conjugate roots of a same irreducible polynomial (for more details see \cite{Perez}).
\end{remark}

\section{The  non-rational case: computing branches and asymptotes}\label{S-new}

Throughout this paper so far, we have dealt with  algebraic plane curves
implicitly and rational parametrically defined. In this section, we present all the previous concepts introduced before for the case of meromorphic functions. In addition, we  present a method
to compute infinity branches and g-asymptotes for these type of functions. This method  is based on the idea presented in Subsection \ref{Sub-improvement}, where we  show how one easily compute the generalized asymptotes  by only determining some simple limits of rational functions constructed from the rational functions defining the input  parametrization.

\para

We remind that the g-asymptotes of an input curve, $\cal C$, are perfect curves computed from the  infinity branches of $\cal C$. That is, once one has the branch
$$B=\{(z, r(z))\in {\Bbb C}^2: \,z\in {\Bbb
C},\,|z|>M\},$$
where (see equality \ref{Eq-inf-branchn}) $$r(z)=mz+a_1z^{1-n_1/n}+\cdots
+a_kz^{1-n_k/n}+a_{k+1}z^{1-N_{k+1}/N}+\cdots,$$ and $0<N_1<\cdots<N_k\leq N<N_{k+1}<\cdots$, $\gcd(N,N_1,\ldots,N_k)=b$, $N=n \cdot b$, $N_j=n_j\cdot b,\,\,j\in\{1,\ldots,k\}$, the asymptote is obtained by considering the  terms
with non negative exponent in the series $r(z)$. We note that we have simplified the non negative
exponents such that $\gcd(n,n_1,\ldots,n_k)=1$, $0<n_1<n_2<\cdots$, and  $n_k\leq n$. We say that $n$ is the degree of $B$ ($\deg(B)$) and  $N$ is  the ramification index of the branch $B$ ($\nu(B)$).

\para

Additionally, we say that the infinity branch $B$ is associated to the  infinity point $P=(1:m:0),\, m\in\mathbb{C}$ and it holds that   $f(z,r(z))=0$ for $z\in {\Bbb
C},\,|z|>M$ when the curve is implicitly defined. If the curve is defined by the rational parametrization, it holds that $$\limit_{t\rightarrow \tau} {\cal P}(t)=\limit_{z\rightarrow \infty} (z, r(z)),$$
where $\tau$ is a value of $t$ for which $\cP$ is not defined and $\cP^*(\tau)=(1:m:0)$. In the case we are dealing in this section, we do not have an implicit equation so we would use this last characterization.

\para

We also remind that  a g-asymptote of $\cal C$ at $B$ is a perfect curve
that approaches $\cal C$ at $B$. So, we need to compute a perfect curve approaching the input curve $\cal C$, then we will get the g-asymptote. For this purpose, we recall that  a curve   defined by a proper parametrization of the form $(t^{n}, a_{n}t^{n}+a_{n-1}t^{n-1}+\cdots+a_0)$ is always perfect (see Section \ref{S-asymptotes}). So, our purpose is to compute a curve of this form approaching the input curve $\cal C$. That is, we impose the condition that
$$\limit_{t\rightarrow \tau} {\cal P}(t)=\limit_{z\rightarrow \infty} (z^n, a_{n}z^{n}+a_{n-1}z^{n-1}+\cdots+a_0),$$
where $\tau$ is a value of $t$ for which $\cP$ is not defined.

\para

Using this idea in fact we see how to determine the  infinity branches
$$B=\{(z, r(z))\in {\Bbb C}^2: \,z\in {\Bbb
C},\,|z|>M\},$$
where  $$r(z)=mz+a_1z^{1-n_1/n}+\cdots
+a_kz^{1-n_k/n}+a_{k+1}z^{1-N_{k+1}/N}+\cdots,$$ and $0<N_1<\cdots<N_k\leq N<N_{k+1}<\cdots$.
Since we can not compute $\varphi(t)$, in order to compute $B$, we use the idea presented in Subsection \ref{Sub-improvement}   and we impose the condition
$$\limit_{t\rightarrow \tau} {\cal P}(t)=\limit_{z\rightarrow \infty} (z, r(z)).$$

In the following, we deal with   meromorphic functions from an open subset of the complex plane onto ${\Bbb C}^n,\,n\geq 1$. We remind that  meromorphic functions are  functions on an open subset $D$ of $\Bbb C$ that are holomorphic on all  $D$ except for a set of isolated points, which are poles of the functions. By abuse of notation,   we denote by ${\cal P}(t):=(p_1(t),\ldots,p_n(t))$ these meromorphic functions and we say that the image of $\cal P$ is a {\it special curve $\cal C$ parametrically defined in   $n$-dimensional space}.

\para

 We develop the method for $n=2$ but it can be straightforward adapted  for dealing with  curves in   $n$-dimensional space. Thus, in this section we have a {\it curve } $\cal C$ that is the image of the  {\it parametrization}
\[{\cal P}(s)=(p_1(s),p_2(s)),\quad p_i(s)=p_{i1}(s)/p_{i2}(s),\quad \gcd(p_{i1},p_{i2})=1,\,\,i=1,2,\]
where the functions $p_{i}(s),\,i\in \{1,2\}$ are meromorphic   on the complex plane.

\para

The different roots of the denominators $p_{i2},\,\i=1,2$ (the poles of the meromorphic functions $p_i(s),\,i=1,2$) yield the infinity branches  associated to the infinity points  $P=(1:m:0),\, m\in\mathbb{C}$. The following theorem shows how to compute the branches, by just computing some simple limits of some  functions constructed from  $\cP(s)$. Afterwards, from each of these branches  we can easily  obtain the g-asymptote.
\para

We assume that   we have  prepared the input curve ${\cal C}$, by means
of a suitable   projective   change of coordinates (if necessary) such that
$(0:1:0)$ is not a point  at infinity. 

\para

\begin{theorem}\label{T-final-ramas-G} Let  $\cal C$ be a curve defined  by a parametrization
\[{\cal P}(s)=(p_1(s),p_2(s)),\quad p_i(s)=p_{i1}(s)/p_{i2}(s),\quad \gcd(p_{i1},p_{i2})=1,\,\,i=1,2.\]
Let $\tau\in {\Bbb C}$ be such that  $p_{i2}(t)=(t-\tau)^{n_i/m_i}\overline{p}_{i2}(t)$ where $\overline{p}_{i2}(\tau)\not=0,\,i=1,2$, and $n_1/m_1\geq 1$. Let us assume that $p_{i1}(t)=(t-\tau)^{u_i/v_i}\overline{p}_{i1}(t),\,i=1,2$ and $0\leq u_1/v_1<n_1/m_1$ and $0\leq u_2/v_2\leq n_2/m_2$.  Let $\gamma:=\lcm(m_1,m_2,v_1,v_2)$ and  $${\cal M}(s)={\cal P}(s^\gamma)=(\wp_1(s),\wp_2(s)),\quad \wp_i(s)=\wp_{i1}(s)/\wp_{i2}(s),\,\,i=1,2.$$ Let $\bar{n}_i:=n_i\gamma/m_i-u_i\gamma/v_i,\,i=1,2$.  An
infinity branch associated to an infinity point $P=(1:m:0), m \in\mathbb{C}$ is given as
$$B=\{(z,r(z))\in
{\Bbb C}^2: \,z\in {\Bbb C},\,|z|>M\},$$
where
\[r(z)=a_{\bar{n}_2}z^{\bar{n}_2/\bar{n}_1}+a_{\bar{n}_2-1}z^{(\bar{n}_2-1)/\bar{n}_1}+\ldots+a_{0}+a_{-1}z^{-1/\bar{n}_1}+a_{-2}z^{-2/\bar{n}_1}+\cdots\]
($m=a_1$), and
$$\begin{array}{ll}
   a_{\bar{n}_2}=\limit_{t\rightarrow \tau}\,\, \displaystyle\frac{\wp_2(t)}{\wp_{1}(t)^{\bar{n}_2/\bar{n}_1}}& \\
  a_{\bar{n}_2-1}=\limit_{t\rightarrow \tau}\, \wp_1(t)^{1/\bar{n}_1}f_1(t),\quad &f_1(t):=\displaystyle \frac{\wp_2(t)}{\wp_{1}(t)^{\bar{n}_2/\bar{n}_1}}-a_{\bar{n}_2}  \\
 a_{\bar{n}_2-2}=\limit_{t\rightarrow \tau}\, \wp_1(t)^{1/\bar{n}_1}f_2(t),\quad &f_2(t):= \wp_1(t)^{1/\bar{n}_1}f_1(t)-a_{\bar{n}_2-1}  \\
\qquad\quad \vdots & \qquad\quad \vdots\\
   a_{\bar{n}_2-i}=\limit_{t\rightarrow \tau}\, \wp_1(t)^{1/\bar{n}_1}f_i(t),\quad &f_i(t):= \wp_1(t)^{1/\bar{n}_1}f_{i-1}(t)-a_{\bar{n}_2-(i-1)},\,\, \\
 \end{array}$$
for $i\in\{2,\ldots,\bar{n}_2,\bar{n}_2+1,\bar{n}_2+2,\cdots\}$.
\end{theorem}
\begin{proof} We use the equality
$$\limit_{t\rightarrow \tau} {\cal P}(t)=\limit_{z\rightarrow \infty} (z,r(z)) =\limit_{z\rightarrow \infty} (z^{\bar{n}_1},r(z^{\bar{n}_1})) $$
  for easily obtaining the coefficients, $a_{n_2-i}$ for $i\in\{2,\ldots,\bar{n}_2,\bar{n}_2+1,\bar{n}_2+2,\cdots\}$. Indeed, we have that
\[\limit_{t\rightarrow \tau}\,\frac{\wp_2(t)}{\wp_{1}(t)^{\bar{n}_2/\bar{n}_1}}=\]\[\limit_{z\rightarrow \infty}\frac{a_{\bar{n}_2}z^{\bar{n}_2/\bar{n}_1}+a_{\bar{n}_2-1}z^{(\bar{n}_2-1)/\bar{n}_1}+\ldots+a_{0}+a_{-1}z^{-1/\bar{n}_1}+a_{-2}z^{-2/\bar{n}_1}+\cdots}{(z^{\bar{n}_1})^{\bar{n}_2/\bar{n}_1}}=a_{\bar{n}_2},\]
\[\limit_{t\rightarrow \tau}\,\frac{\wp_2(t)-a_{\bar{n}_2}\wp_{1}(t)^{\bar{n}_2/\bar{n}_1}}{\wp_{1}(t)^{(\bar{n}_2-1)/\bar{n}_1}}=\]
\[\limit_{t\rightarrow \infty}\frac{a_{\bar{n}_2-1}z^{(\bar{n}_2-1)/\bar{n}_1}+\ldots+a_{0}+a_{-1}z^{-1/\bar{n}_1}+a_{-2}z^{-2/\bar{n}_1}+\cdots}{(z^{\bar{n}_1})^{(\bar{n}_2-1)/\bar{n}_1}}=a_{\bar{n}_2-1},\]
\[\limit_{t\rightarrow \tau}\,\frac{\wp_2(t)-a_{\bar{n}_2}\wp_{1}(t)^{\bar{n}_2/\bar{n}_1}-a_{\bar{n}_2-1}\wp_{1}(t)^{(\bar{n}_2-1)/\bar{n}_1}}{\wp_{1}(t)^{(\bar{n}_2-2)/\bar{n}_1}}=\]\[\limit_{z\rightarrow \infty}\frac{a_{\bar{n}_2-2}z^{(\bar{n}_2-2)/\bar{n}_1}+\ldots+a_{0}+a_{-1}z^{-1/\bar{n}_1}+a_{-2}z^{-2/\bar{n}_1}+\cdots}{(z^{\bar{n}_1})^{(\bar{n}_2-2)/n_1}}=a_{\bar{n}_2-2}\]
and, reasoning similarly,
\[\limit_{t\rightarrow \tau}\,(\wp_2(t)-a_{\bar{n}_2}\wp_{1}(t)^{n_2/n_1}-a_{\bar{n}_2-1}\wp_{1}(t)^{(\bar{n}_2-1)/\bar{n}_1}-\cdots-a_1\wp_{1}(t)^{1/\bar{n}_1})=\]\[\limit_{z \rightarrow \infty} (a_{0}+a_{-1}z^{-1/\bar{n}_1}+a_{-2}z^{-2/\bar{n}_1}+\cdots)=a_{0}.\]
Additionally, we have that
\[\limit_{t\rightarrow \tau}\,\frac{\wp_2(t)-a_{\bar{n}_2}\wp_{1}(t)^{n_2/n_1}-a_{\bar{n}_2-1}\wp_{1}(t)^{(\bar{n}_2-1)/\bar{n}_1}-\cdots-a_0}{\wp_1(t)^{-1/n_1}}=\]\[\limit_{z \rightarrow \infty}\frac{a_{-1}z^{-1/\bar{n}_1}+a_{-2}z^{-2/\bar{n}_1}+\cdots}{(z^{n_1})^{-1/n_1}}=a_{-1},\]
\[\limit_{t\rightarrow \tau}\,\frac{\wp_2(t)-a_{\bar{n}_2}\wp_{1}(t)^{n_2/n_1}-a_{\bar{n}_2-1}\wp_{1}(t)^{(\bar{n}_2-1)/\bar{n}_1}-\cdots-a_0-a_{-1}z^{-1/\bar{n}_1}}{\wp_1(t)^{-2/n_1}}=\]\[\limit_{z \rightarrow \infty}\frac{a_{-2}z^{-2/\bar{n}_1}+\cdots}{(z^{\bar{n}_1})^{-2/\bar{n}_1}}=a_{-2},\]
and, reasoning similarly for $i\geq 3$
\[\limit_{t\rightarrow \tau}\,\frac{\wp_2(t)-a_{\bar{n}_2}\wp_{1}(t)^{n_2/n_1}-a_{\bar{n}_2-1}\wp_{1}(t)^{(\bar{n}_2-1)/\bar{n}_1}-\cdots-a_0-a_{-1}z^{-1/\bar{n}_1}-\cdots-a_{-(i-1)}z^{-(i-1)/\bar{n}_1}}{\wp_1(t)^{-i/\bar{n}_1}}=\]\[\limit_{z \rightarrow \infty}\frac{a_{-i}z^{-i/\bar{n}_1}+\cdots}{(z^{n_1})^{-i/n_1}}=a_{-i}.\]
Finally, we observe that we may write
 \[\frac{\wp_2(t)-a_{\bar{n}_2}\wp_{1}(t)^{\bar{n}_2/\bar{n}_1}}{\wp_{1}(t)^{(\bar{n}_2-1)/\bar{n}_1}}=\wp_1(t)^{1/\bar{n}_1}f_1(t),\quad f_1(t):=\displaystyle \frac{\wp_2(t)}{\wp_{1}(t)^{\bar{n}_2/\bar{n}_1}}-a_{\bar{n}_2}\]
\[\frac{\wp_2(t)-a_{\bar{n}_2}p_{1}(t)^{\bar{n}_2/\bar{n}_1}-a_{\bar{n}_2-1}p_{1}(t)^{(\bar{n}_2-1)/\bar{n}_1}}{\wp_{1}(t)^{(\bar{n}_2-2)/\bar{n}_1}}=\wp_1(t)^{1/\bar{n}_1}f_2(t),\quad f_2(t):= \wp_1(t)^{1/\bar{n}_1}f_1(t)-a_{\bar{n}_2-1}\]
and in general,
$$\frac{\wp_2(t)-a_{\bar{n}_2}\wp_{1}(t)^{n_2/n_1}-a_{\bar{n}_2-1}\wp_{1}(t)^{(\bar{n}_2-1)/\bar{n}_1}-\cdots-a_0-a_{-1}z^{-1/\bar{n}_1}-\cdots-a_{(\bar{n}_2-i+1)}z^{(\bar{n}_2-i+1)/\bar{n}_1}}{\wp_1(t)^{(\bar{n}_2-i)/\bar{n}_1}}$$$$=\wp_1(t)^{1/\bar{n}_1}f_i(t),\quad  f_i(t):= \wp_1(t)^{1/\bar{n}_1}f_{i-1}(t)-a_{\bar{n}_2-i+1}$$
for $i\in\{2,\ldots,\bar{n}_2,\bar{n}_2+1,\bar{n}_2+2,\cdots\}$.
\end{proof}

\para

\begin{remark}\label{R-degree}
Note that $\bar{n}_1$ is not the degree of $B$ but its ramification index. In addition, we observe  that since $(0:1:0)$ is not a point  at infinity then $\bar{n}_1\geq \bar{n}_2$.
\end{remark}

\para

\begin{remark}\label{R-infpoint2-G}
From the above construction, each root $\tau$ of $p_{12}(t)$ yields an infinity branch and, hence, an infinity point $P=(1:m:0)$ (see Remark \ref{R-infpoint}). Note that the parametrization $\mathcal{M}(t)$ constructed in the previous theorem can be expressed as
$\mathcal{M}(t)=\left(\frac{\wp_{11}(t)}{\wp(t)},\frac{\wp_{12}(t)}{\wp(t)}\right)$,
where $\wp(t)=\lcm(\wp_{12}(t),\wp_{22}(t))$ and $\wp_{1i}(t)=\wp_i(t)\wp(t)$. Now, the corresponding projective curve is parametrized by $\mathcal{M}^*(t)=(\wp_{11}(t),\wp_{12}(t),\wp(t))$
and the infinity point associated to $\tau$ is $P=(\wp_{11}(\tau):\wp_{12}(\tau):0)$. Note that $P=(1:\wp_{12}(\tau)/\wp_{11}(\tau):0)$ if $\wp_{11}(\tau)\not=0$.
\end{remark}

\para

From Theorem  \ref{T-final-ramas-G}, we easily get the following theorem that allows us to compute the g-asymptote.

\para

\begin{theorem}\label{T-final-G} Let  $\cal C$ be a curve defined  by a parametrization
\[{\cal P}(s)=(p_1(s),p_2(s)),\quad p_i(s)=p_{i1}(s)/p_{i2}(s),\quad \gcd(p_{i1},p_{i2})=1,\,\,i=1,2.\]
Let $\tau\in {\Bbb C}$ be such that  $p_{i2}(t)=(t-\tau)^{n_i/m_i}\overline{p}_{i2}(t)$ where $\overline{p}_{i2}(\tau)\not=0,\,i=1,2$, and $n_1/m_1\geq 1$, and let $B$ be the corresponding infinity branch. Let us assume that $p_{i1}(t)=(t-\tau)^{u_i/v_i}\overline{p}_{i1}(t),\,i=1,2$ and $0\leq u_1/v_1<n_1/m_1$ and $0\leq u_2/v_2\leq n_2/m_2$.  Let $\gamma:=\lcm(m_1,m_2,v_1,v_2)$ and  $${\cal M}(s)={\cal P}(s^\gamma)=(\wp_1(s),\wp_2(s)),\quad \wp_i(s)=\wp_{i1}(s)/\wp_{i2}(s),\,\,i=1,2.$$ Let $\bar{n}_i:=n_i\gamma/m_i-u_i\gamma/v_i,\,i=1,2$.
 A g-asymptote of an infinity branch $B$ associated to an infinity point $P=(1:m:0), m \in\mathbb{C}$ is defined by a proper reparametrization of
$$\widetilde{\cal Q}(t)=(t^{\bar{n}_1},\,  a_{\bar{n}_2}t^{\bar{n}_2}+a_{\bar{n}_2-1}t^{\bar{n}_2-1}+\ldots+a_{0}),$$
($m=a_1$), where
$$\begin{array}{ll}
   a_{\bar{n}_2}=\limit_{t\rightarrow \tau}\,\, \displaystyle\frac{\wp_2(t)}{\wp_{1}(t)^{\bar{n}_2/\bar{n}_1}}& \\
  a_{\bar{n}_2-1}=\limit_{t\rightarrow \tau}\, \wp_1(t)^{1/\bar{n}_1}f_1(t),\quad &f_1(t):=\displaystyle \frac{\wp_2(t)}{\wp_{1}(t)^{\bar{n}_2/\bar{n}_1}}-a_{\bar{n}_2}  \\
 a_{\bar{n}_2-2}=\limit_{t\rightarrow \tau}\, \wp_1(t)^{1/\bar{n}_1}f_2(t),\quad &f_2(t):= \wp_1(t)^{1/\bar{n}_1}f_1(t)-a_{\bar{n}_2-1}  \\
\qquad\quad \vdots & \qquad\quad \vdots\\
   a_{\bar{n}_2-i}=\limit_{t\rightarrow \tau}\, \wp_1(t)^{1/\bar{n}_1}f_i(t),\quad &f_i(t):= \wp_1(t)^{1/\bar{n}_1}f_{i-1}(t)-a_{\bar{n}_2-(i-1)},\,\,i\in\{2,\ldots,\bar{n}_2\}. \\
 \end{array}$$
\end{theorem}
\begin{proof} We first note that a
  curve $\widetilde{\cal C}$ approaching the infinity branch $B$ is given by
$$\widetilde{Q}(t)=(t^{\bar{n}_1},\widetilde{r}(t^{\bar{n}_1}))=(t^{\bar{n}_1},\,  a_{\bar{n}_2}t^{\bar{n}_2}+a_{\bar{n}_2-1}t^{\bar{n}_2-1}+\ldots+a_{0}),$$
where $\bar{n}_1\geq \bar{n}_2$ (see Remark \ref{R-degree}), and $\widetilde{r}(z)$ can be computed from $r(z)$ by removing the terms with negative exponent (see  Theorem  \ref{T-final-ramas-G} and Section \ref{Sub-asintotasconstructions}). In order to prove that $\widetilde{\cal C}$ is a g-asymptote, we need to prove that it is a perfect curve. For this purpose, we observe that the degree of $\widetilde{\cal C}$ is not necessary $\bar{n}_1$ since $\widetilde{Q}$ could be improper which is equivalent to $\gcd(\bar{n}_1, \bar{n}_2,\ldots, \bar{n}_2-j)\not=1$ for every $j=0,\ldots, \bar{n}_2-1$ such that $a_{\bar{n}_2-j}\not=0$ (see Remark  \ref{R-degree}).
  Let us assume that $\gcd(\bar{n}_1, \bar{n}_2,\ldots,\bar{n}_2-j)=\beta$ for every $j=0,\ldots, \bar{n}_2-1$ such that $a_{\bar{n}_2-j}\not=0$. Then, let we consider the reparametrizatin
  \[\cM(t)= \cP(t^{1/\beta})=(t^{\bar{n}_1/\beta},\,  a_{\bar{n}_2}t^{\bar{n}_2/\beta}+a_{\bar{n}_2-1}t^{(\bar{n}_2-1)/\beta}+\ldots+a_{0})\]
is a proper reparametrization of $\widetilde{Q}$ (note that the exponents are natural integers since  $\gcd(\bar{n}_1, \bar{n}_2,\ldots,\bar{n}_2-j)=\beta$ for every $j=0,\ldots, \bar{n}_2-1$ such that $a_{\bar{n}_2-j}\not=0$). Then we get that the theorem outputs an asymptote since $\widetilde{\cal C}$ is perfect (it has an only branch and the degree of the curve $\widetilde{\cal C}$ which is  ${\bar{n}_1/\beta}$, see Remark  \ref{R-degree}, is equal to the degree of the branch).

\end{proof}

\para

In the following corollaries, we analyze the special case of the vertical and horizontal g-asymptotes, i.e. lines of the form $x-a$ or $y-b$, where $a,b\in {\Bbb C}$ (observe that these asymptotes correspond to branches associated to the infinity points $(0:1:0)$ and $(1:0:0)$, respectively). More precisely, we prove that these asymptotes are obtained from the non--common roots of the denominators of the given parametrization. Note that  in the practical design  of engineering and modeling applications, the   curves  are usually presented by numerical coefficients and $\cP(s)$    mostly satisfies that $\gcd(p_{12}, p_{22})=1$.

\para

\begin{corollary}\label{C-final-G}  Let  $\cal C$ be a curve defined  by a parametrization
\[{\cal P}(s)=(p_1(s),p_2(s),\quad p_i(s)=p_{i1}(s)/p_{i2}(s),\quad \gcd(p_{i1},p_{i2})=1,\,\,i=1,2.\]
Let $\tau\in {\Bbb C}$ be such that  $p_{i2}(t)=(t-\tau)^{n_i/m_i}\overline{p}_{i2}(t)$ where $\overline{p}_{i2}(\tau)\not=0,\,i=1,2$, and $n_1/m_1\geq 1$, and let $B$ be the corresponding infinity branch. Let us assume that $p_{i1}(t)=(t-\tau)^{u_i/v_i}\overline{p}_{i1}(t),\,i=1,2$ and $0\leq u_1/v_1<n_1/m_1$ and $0\leq u_2/v_2=n_2/m_2$.  Let $\gamma:=\lcm(m_1,m_2,v_1,v_2)$ and  $${\cal M}(s)={\cal P}(s^\gamma)=(\wp_1(s),\wp_2(s)),\quad \wp_i(s)=\wp_{i1}(s)/\wp_{i2}(s),\,\,i=1,2.$$ Let $\bar{n}_i:=n_i\gamma/m_i-u_i\gamma/v_i,\,i=1,2$.
It holds that a g-asymptote of $\cal C$ corresponding to the infinity point $(1:0:0)$ is the horizontal line  $y-\wp_{2}(\tau)=0$, defined by the parametrization
$\widetilde{\cal Q}(t)=(t,\,  \wp_{2}(\tau)).$
   \end{corollary}
   \begin{proof}
     We apply Theorem \ref{T-final} with $\bar{n}_2=0.$
   \end{proof}

\para

\begin{corollary}\label{C-final2-G}   Let  $\cal C$ be a curve defined  by a parametrization
\[{\cal P}(s)=(p_1(s),p_2(s),\quad p_i(s)=p_{i1}(s)/p_{i2}(s),\quad \gcd(p_{i1},p_{i2})=1,\,\,i=1,2.\]
Let $\tau\in {\Bbb C}$ be such that  $p_{i2}(t)=(t-\tau)^{n_i/m_i}\overline{p}_{i2}(t)$ where $\overline{p}_{i2}(\tau)\not=0,\,i=1,2$, and $n_2/m_2\geq 1$, and let $B$ be the corresponding infinity branch. Let us assume that $p_{i1}(t)=(t-\tau)^{u_i/v_i}\overline{p}_{i1}(t),\,i=1,2$ and $0\leq u_2/v_2<n_2/m_2$ and $0\leq u_1/v_1=n_1/m_1$.  Let $\gamma:=\lcm(m_1,m_2,v_1,v_2)$ and  $${\cal M}(s)={\cal P}(s^\gamma)=(\wp_1(s),\wp_2(s)),\quad \wp_i(s)=\wp_{i1}(s)/\wp_{i2}(s),\,\,i=1,2.$$ Let $\bar{n}_i:=n_i\gamma/m_i-u_i\gamma/v_i,\,i=1,2$.
It holds that a g-asymptote of $\cal C$ corresponding to the infinity point $(0:1:0)$  is the vertical line  $x-p_{1}(\tau)=0$, defined by the parametrization
$\widetilde{\cal Q}(t)=(\wp_{1}(\tau),\,  t).$
   \end{corollary}
   \begin{proof} We apply Corollary \ref{C-final-G} to the parametrization $(p_2(s), p_1(s))$ and we get the asymptote defined by the parametrization $(t, p_{1}(\tau))$. Afterwards, we undo  the change of coordinates (see statement 2 in Remark \ref{R-param-l}).
   \end{proof}
   \para




In the following we introduce the Algorithm {\sf  Asymptotes Construction-Parametric Non-Rational Case}, which uses the above results for computing the g-asymptotes of a  plane curve.

\begin{center}
\fbox{\hspace*{2 mm}\parbox{15cm}{ \vspace*{2 mm} {\bf Algorithm
{\sf    Asymptotes Construction-Parametric Non-Rational Case.}}
\vspace*{0.2cm}

\noindent {\sf Given} a   curve
$\cal C$ defined by
${\cal P}(s)=(p_1(s),p_2(s)),\quad p_i(s)=p_{i1}(s)/p_{i2}(s),\quad \gcd(p_{i1},p_{i2})=1,\,\,i=1,2$, where $p_{i}(s),\,i\in \{1,2\}$ are  meromorphic functions on the complex plane,
  the algorithm {\sf outputs}  one asymptote
for each of its infinity branches.
\begin{itemize}

\item[1.]  Let  $\tau_1,\ldots, \tau_k\in {\Bbb C}$ be the roots of $p_{12}$  that is, the poles of $p_1$.

\item[2.]For each   $\tau_i$, $i\in\{1,\ldots,k\},$ do:
\begin{itemize}

\item[2.1] Write  $p_{i2}(t)=(t-\tau)^{n_i/m_i}\overline{p}_{i2}(t)$ where $\overline{p}_{i2}(\tau)\not=0,\,i=1,2$, and $n_1/m_1\geq 1$, and $p_{i1}(t)=(t-\tau)^{u_i/v_i}\overline{p}_{i1}(t),\,i=1,2$.
Check whether  $0\leq u_1/v_1<n_1/m_1$ and $0\leq u_2/v_2\leq n_2/m_2$ and in the affirmative case consider $\gamma:=\lcm(m_1,m_2,v_1,v_2)$ and  $${\cal M}(s)={\cal P}(s^\gamma)=(\wp_1(s),\wp_2(s)),\quad \wp_i(s)=\wp_{i1}(s)/\wp_{i2}(s),\,\,i=1,2.$$ Let $\bar{n}_i:=n_i\gamma/m_i-u_i\gamma/v_i,\,i=1,2$.

\item[2.2.] Compute
$$\begin{array}{ll}
   a_{\bar{n}_{2i}}=\limit_{t\rightarrow \tau_i}\,\, \frac{\wp_2(t)}{p_{1}(t)^{\bar{n}_{2i}/\bar{n}_{1i}}}& \\
  a_{\bar{n}_{2i}-1}=\limit_{t\rightarrow \tau_i}\, \wp_1(t)^{1/n_{1i}}f_1(t),\quad &f_1(t):= \frac{\wp_2(t)}{\wp_{1}(t)^{\bar{n}_{2i}/\bar{n}_{1i}}}-a_{\bar{n}_{2i}}  \\
 a_{\bar{n}_{2i}-2}=\limit_{t\rightarrow \tau_i}\, \wp_1(t)^{1/n_{1i}}f_2(t),\quad &f_2(t):= \wp_1(t)^{1/\bar{n}_{1i}}f_1(t)-a_{\bar{n}_{2i}-1}  \\
 \qquad \quad \vdots & \qquad \quad \vdots\\
   a_{\bar{n}_{2i}-j}=\limit_{t\rightarrow \tau_i}\, \wp_1(t)^{1/n_{1i}}f_j(t),\quad &f_j(t):= \wp_1(t)^{1/\bar{n}_{1i}}f_{j-1}(t)-a_{\bar{n}_{2i}-(j-1)}, \\
 \end{array}$$
for $j\in\{2,\ldots, \bar{n}_{2i}\}$.

\item[2.3.] Let $\widetilde{\cal C}_i$ be the asymptote  defined by the proper parametrization
$$\widetilde{\cal Q}_i(t)=(t^{\bar{n}_{1i}},\,  a_{\bar{n}_{2i}}t^{\bar{n}_{2i}}+a_{\bar{n}_{2i}-1}t^{\bar{n}_{2i}-1}+\ldots+a_{0})
 \in
{\Bbb C}[t]^2.$$
\end{itemize}

\item[3.]  If there exist  $s_1,\ldots, s_l\in {\Bbb C}$   roots of $p_{22}(s)$ (that is, the poles of $p_2$) such that $p_{12}(s_j)\not=0$ for $j\in\{1,\ldots, l\}$ then let
$\widetilde{\cal D}_i$ be  the vertical asymptote defined by the proper parametrization
$$\widetilde{\cal Q}_i(t)=(\wp_1(s_i), t) \in
{\Bbb C}[t]^2,\,\,i\in\{1,\ldots,l\}.$$

\item[4.]  {\sf Return} the asymptotes $\widetilde{\cal C}_1,\ldots,\widetilde{\cal C}_k$ and $\widetilde{\cal D}_1,\ldots,\widetilde{\cal D}_l$.

\end{itemize} }\hspace{2 mm}}
\end{center}

\para

By applying Algorithm {\sf  Asymptotes Construction-Parametric Non-Rational Case}, we can easily obtain all the g-asymptotes of any   plane curve, as  the following examples show.

\para

\begin{example}\label{E-method-ant-b-G} We consider the curve ${\cal C}$   defined by the
parametrization
$${\cal P}(s)=(p_1(s), p_2(s))=\left(\frac{\,s^{1/2}+1}{\,s^{1/2} \sin(s)}, \frac{s^2+s+5}{\sin(s)}\right).$$ We apply the algorithm  {\sf  Asymptotes Construction-Parametric Non-Rational Case.}\para

\para

\noindent{\sf Step 1:}  We observe that $p_{12}(s)$ and  $p_{22}(s)$ have the root  $\tau=0$, and $p_{i}(s)=p_{i1}/p_{i2},\,i\in \{1,2\}$ are  meromorphic functions on the complex plane.

\para

\noindent{\sf Step 2:} For $\tau=0:$ \\

\noindent{\sf Step 2.1:} We  write
$$p_{11}(s)=s^{1/2}+1,\quad p_{12}(s)=s^{3/2}(1-1/6 s^{2}+1/120 s^4-1/5040 s^6+\cdots)$$
$$p_{21}(s)=s^2+s+5,\quad p_{22}(s)=s- 1/6 s^3+ 1/120 s^5- 1/5040 s^7+\cdots.$$
Observe that $p_{i1}(0)\not=0$ for $i=1,2$ and thus, we consider
 $${\cal M}(s)={\cal P}(s^2)=(\wp_1(s),\wp_2(s))=\left(\frac{s+1}{s \sin(s^2)}, \frac{s^4+s^2+5}{\sin(s^2)}\right),\quad \wp_i(s)=\wp_{i1}(s)/\wp_{i2}(s),\,\,i=1,2.$$
We get that $\bar{n}_1:=3$ and $\bar{n}_2:=2$.\\

\noindent{\sf Step 2.2:} We have that
$$\begin{array}{ll}
   a_{2}=\limit_{t\rightarrow 2}\,\frac{p_2(t)}{p_{1}(t)^{2}}=5 &\\
  a_{1}=\limit_{t\rightarrow 2}\, p_1(t)f_1(t)-10/3,\quad &f_1(t):=\frac{p_2(t)}{p_{1}(t)^{2}}-a_{2} \\
 a_{0}=\limit_{t\rightarrow 2}\, p_1(t) f_2(t)=8/3,\quad &f_2(t):=p_1(t)f_{1}(t)-a_{1}. \\
 \end{array}$$

\noindent{\sf Step 2.3:} We  obtain the asymptote $\widetilde{\cal C}$, defined by the proper parametrization
$$\widetilde{\cal Q}(t)=(t^3,\,  5t^2-10/3t+8/3).$$

\para

\noindent{\sf Step 3:} We observe that all the roots of $p_{12}(s)$ are also roots of $p_{22}(s)$ (that is, the poles of $p_1$ are the poles of $p_2$), so there are no vertical asymptotes.

\para

\noindent{\sf Step 4:} The algorithm {\sf returns} the   asymptote  $\widetilde{\cal C}$ of the input curve, $\cal C$ (see Figure \ref{F-ex1-method-ant-b-G}).

\begin{figure}[h]
$$
\begin{array}{cc}
\psfig{figure=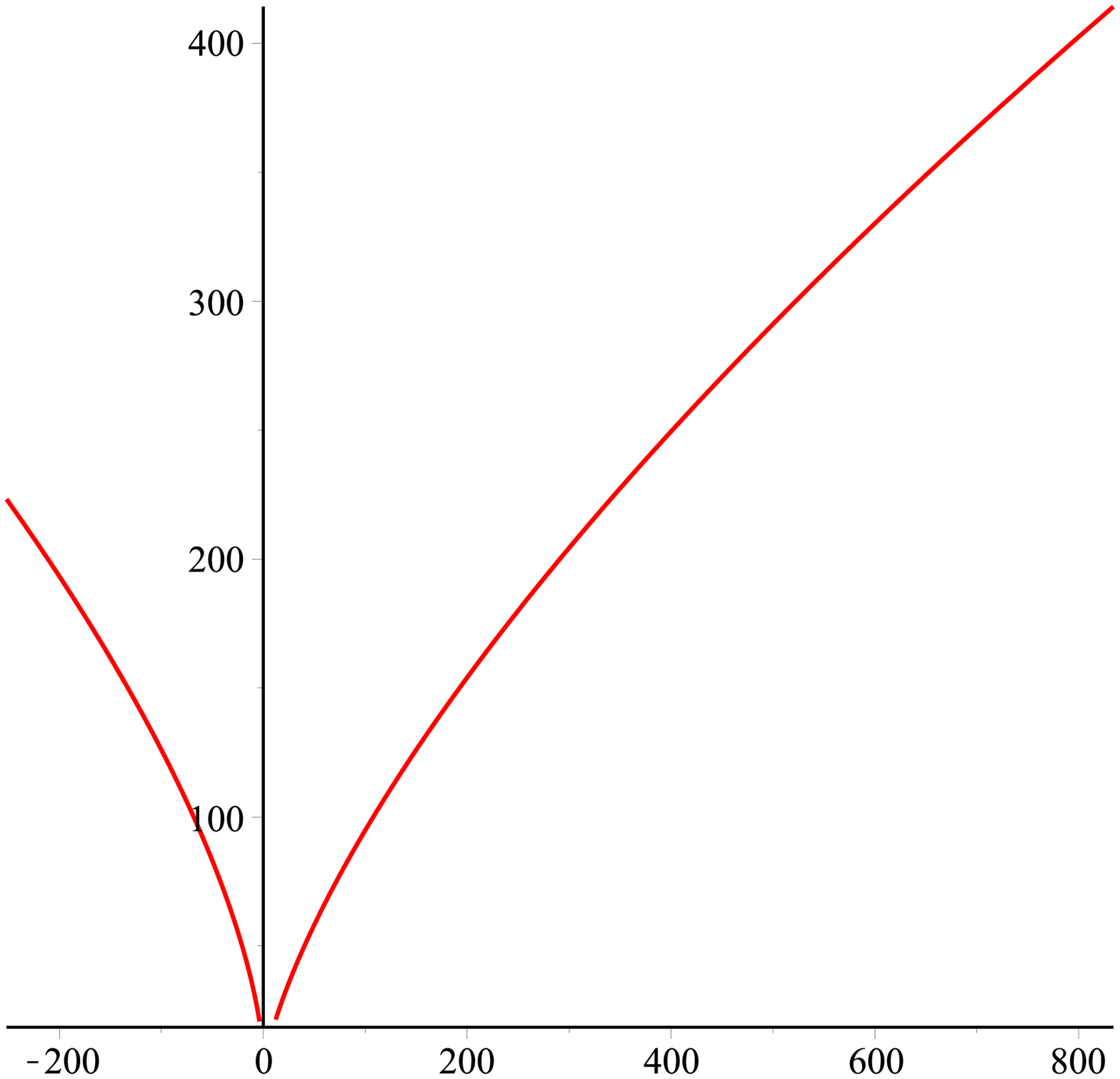,width=6 cm,height=6 cm} &
\psfig{figure=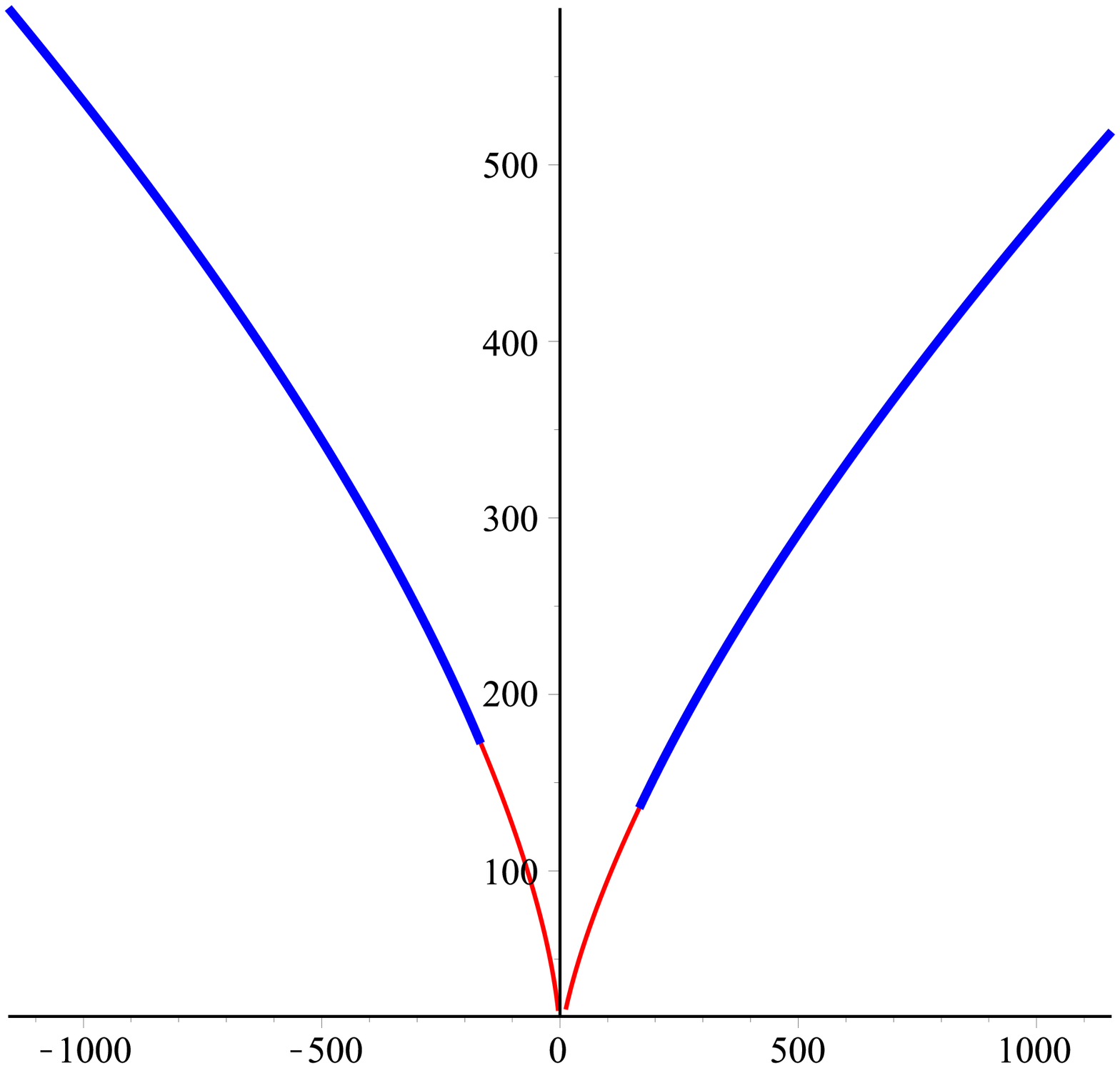,width=6 cm,height=6 cm}
\end{array}
$$ \caption{Curve $\cal C$ (left) and curve and asymptotes (right).}\label{F-ex1-method-ant-b-G}
\end{figure}

\end{example}

\para

\para

\begin{example}\label{E-method-ant-b-G-2} We consider the curve ${\cal C}$   defined by the
parametrization
$${\cal P}(s)=(p_1(s), p_2(s))=\left(\frac{\cos(s)-1}{s^{3/2}\sin(s)}, \frac{\sin(s\pi)}{s^{1/2} \sin(s)}\right).$$ We apply the algorithm  {\sf  Asymptotes Construction-Parametric Non-Rational Case.}\para

\para

\noindent{\sf Step 1:}  We observe that $p_{12}(s)$ and  $p_{22}(s)$ have the root  $\tau=0$, and $p_{i}(s)=p_{i1}/p_{i2},\,i\in \{1,2\}$ are  meromorphic functions on the complex plane.

\para

\noindent{\sf Step 2:} For $\tau=0:$ \\

\noindent{\sf Step 2.1:} We  write
$$p_{12}(s)=s^{5/2}(1-1/6 s^2+1/120s^4-1/5040 s^6+\cdots)$$
$$p_{22}(s)=s^{3/2}(1-1/6 s^2+1/120 s^4-1/5040s^6+\cdots)$$
$$p_{11}(s)s^2(-1/2+1/24s^2-1/720s^4+1/40320s^6+\cdots)$$
$$p_{21}(s)=\pi s (1-1/6s^2\pi^2+1/120s^4 \pi^4-1/5040 s^5\pi^6+\cdots).$$
Thus, we consider
 $${\cal M}(s)={\cal P}(s^2)=(\wp_1(s),\wp_2(s))=\left(\frac{\cos(s^2)-1}{s^3 \sin(s^2)}, \frac{\sin(s^2 \pi)}{s \sin(s^2)}\right),\quad \wp_i(s)=\wp_{i1}(s)/\wp_{i2}(s),\,\,i=1,2.$$
We get that $\bar{n}_1:=1$ and $\bar{n}_2:=1$.\\

\noindent{\sf Step 2.2:} We have that
$$
   a_{1}=\limit_{t\rightarrow 2}\,\frac{p_2(t)}{p_{1}(t)}=-2\pi,\qquad
  a_{0}=\limit_{t\rightarrow 2}\, p_1(t)f_1(t)=0,\,\, f_1(t):=\frac{p_2(t)}{p_{1}(t)}-a_{1}. \\
 $$

\noindent{\sf Step 2.3:} We  obtain the asymptote $\widetilde{\cal C}$, defined by the proper parametrization
$$\widetilde{\cal Q}(t)=(t,\,  -2\pi t).$$

\para

\noindent{\sf Step 3:} We observe that all the roots of $p_{12}(s)$ are also roots of $p_{22}(s)$ (that is, the  the poles of $p_1$ are the poles of $p_2$), so there are no vertical asymptotes.

\para

\noindent{\sf Step 4:} The algorithm {\sf returns} the   asymptote  $\widetilde{\cal C}$ of the input curve, $\cal C$ (see Figure \ref{F-ex1-method-ant-b-G2}).

\begin{figure}[h]
$$
\begin{array}{cc}
\psfig{figure=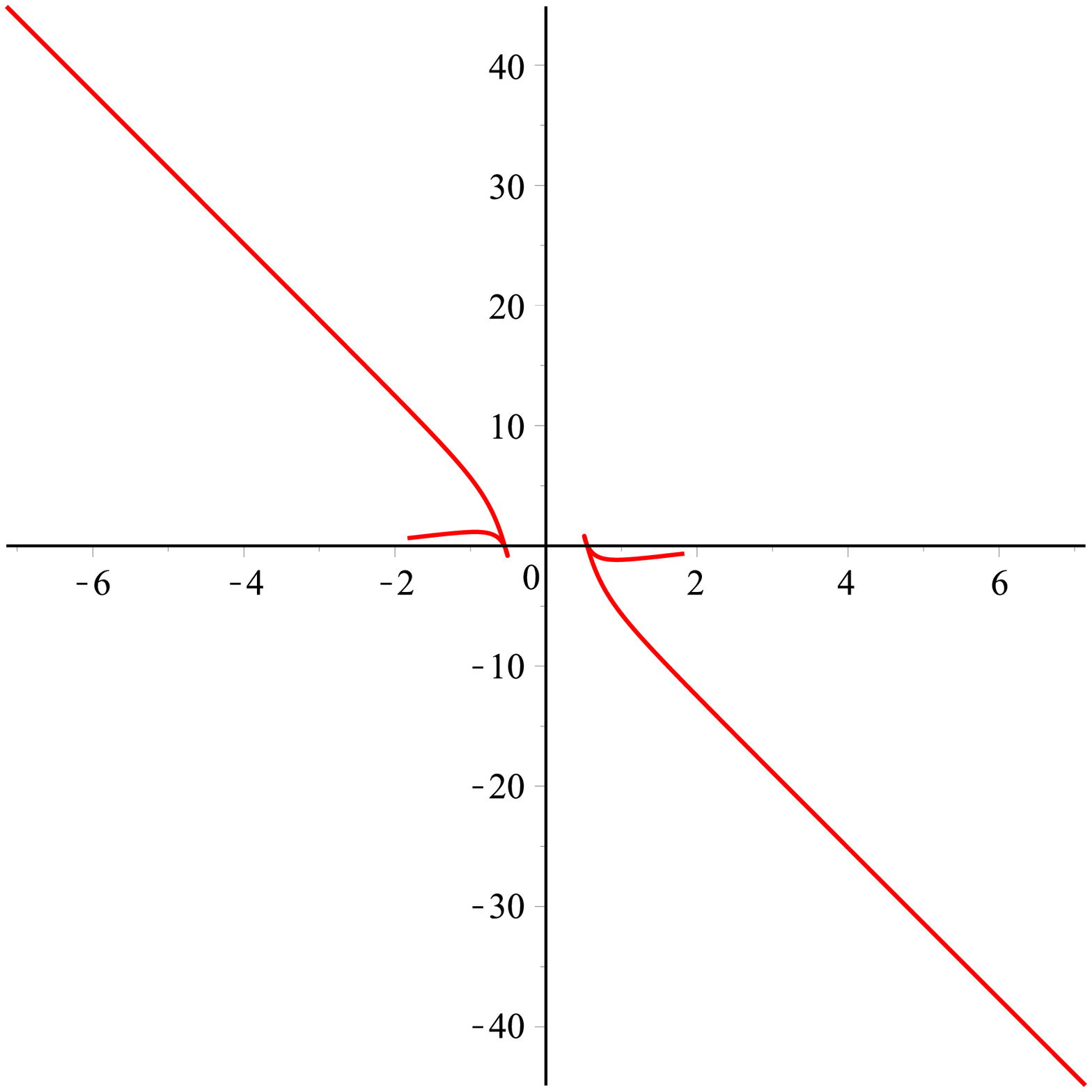,width=6 cm,height=6 cm} &
\psfig{figure=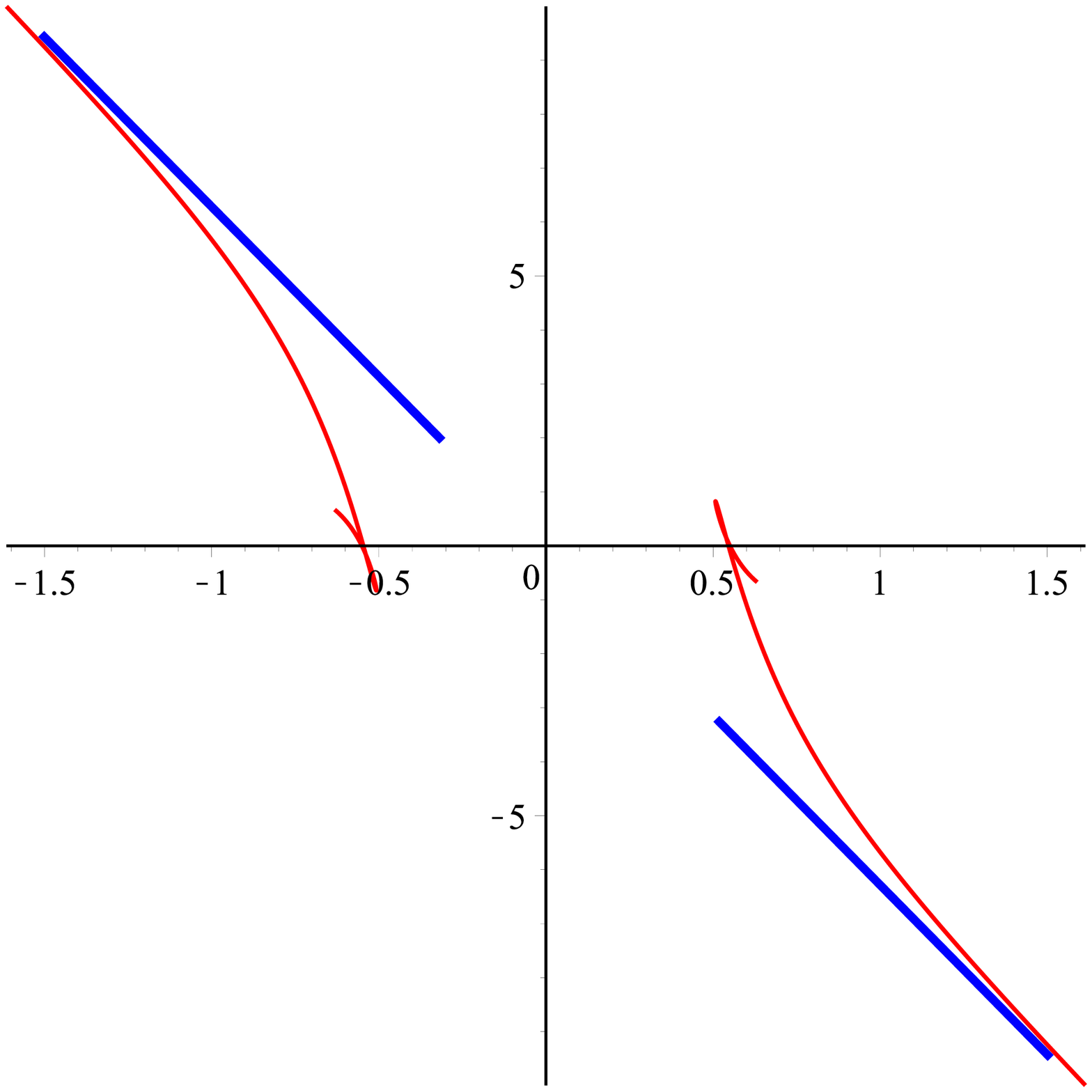,width=6 cm,height=6 cm}
\end{array}
$$ \caption{Curve $\cal C$ (left) and curve and asymptotes (right).}\label{F-ex1-method-ant-b-G2}
\end{figure}

\end{example}

\para

\begin{remark}
The method above described may be straightforward adapted for dealing with  algebraic curves in  the $n-$dimensional space. For instance, if $n=3$, we have a parametrization
${\cal P}(s)=(p_1(s),p_2(s),p_3(s))$ with $p_i(s)=p_{i1}(s)/p_{i2}(s),\,i=1,2,3$ meromorphic functions and $\gcd(p_{i1},p_{i2})=1,\,\,i=1,2,3$. Then, the asymptotes are $$\widetilde{\cal Q}=(t^{n_1},\, a_{n_{2}}t^{n_{2}}+a_{n_{2}-1}t^{n_{2}-1}+\ldots+a_{0},\, b_{m_{2}}t^{m_{2}}+b_{m_{2}-1}t^{m_{2}-1}+\ldots+b_{0}).$$
These asymptotes can be computed by successively applying the algorithm to each component of $\cP$. Note that as in the planar case, roots $\tau\in {\Bbb C}$ such that $p_{22}(\tau)=p_{32}(\tau)=0\neq p_{12}(\tau)$ could appear (see Step 3 of the algorithm). In this case we must look for asymptotes of the form $\widetilde{\cal M}=(m_1(t), t^{n}, m_3(t))$ or $\widetilde{\cal M}=(m_1(t), m_2(t), t^{n})$. Example \ref{ex-method-general-n-dimensional} illustrates these ideas.
\end{remark}

\para

\begin{example}\label{ex-method-general-n-dimensional}
Let ${\cal C}$ be the space curve defined by the
parametrization
$${\cal P}(s)=(p_1(s), p_2(s))=\left(\frac{(\sqrt{s+1}+2)^{1/4}}{s (1+(\sqrt{s+1}+2)^{3/4})},  \frac{\sqrt{\sqrt{s+1}+2}}{ 1+(\sqrt{s+1}+2)^{3/4}},  \frac{s+3}{\sin(s)}\right).$$  We apply the algorithm  {\sf Asymptotes Construction-Parametric Non-Rational Case} for obtaining the different g-asymptotes. \para

\para

\noindent{\sf Step 1:}  We observe that $p_{12}(s)$  and  $p_{32}(s)$  have  the root  $\tau_1=0$  and  $p_{12}(s)$  and  $p_{22}(s)$  have  the root $\tau_2=\alpha$ where $1+(\,(\alpha+1)^{1/2}+2)^{3/4}=0$. Note that  $p_{i}(s)=p_{i1}/p_{i2},\,i\in \{1,2\}$ are  meromorphic functions on the complex plane.

\para

\noindent{\sf Step 2:} \

\begin{itemize}
\item For $\tau_1=0:$ \\

\noindent{\sf Step 2.1:} We  have that $\bar{n}_{11}=1$ and $\bar{n}_{13}=1$  and $\bar{n}_{12}=0$. Note that $p_{32}(s)=s(1-1/6s^2+1/120s^4-1/5040s^6+\cdots)$ and $p_{i1}(0)\not=0$ for $i=1,2,3$.\\

\noindent{\sf Step 2.2:} We have that
$$\begin{array}{ll}
   a_{1}=\limit_{t\rightarrow 2}\,\frac{p_3(t)}{p_{1}(t)}=3^{3/4}+3^{3/2} &\\
  a_{0}=\limit_{t\rightarrow 2}\, p_1(t)f_1(t)=\frac{10\cdot 3^{3/4}+7}{8(1+3^{3/4})},\quad &f_1(t):=\frac{p_3(t)}{p_{1}(t)^{2}}-a_{1}. \\
 \end{array}$$

\noindent{\sf Step 2.3:} We  obtain the asymptote $\widetilde{\cal C}_1$, defined by the proper parametrization
$$\widetilde{\cal Q}_1(t)=\left(t,\, p_2(0),\, (3^{3/4}+3^{3/2})t+\frac{10\cdot 3^{3/4}+7}{8(1+3^{3/4})}\right)$$$$=\left(t,\,\frac{\,3^{1/2}}{1+3^{3/4}},\, (3^{3/4}+3^{3/2})t+\frac{10\cdot 3^{3/4}+7}{8(1+3^{3/4})}\right).$$

\item For  $\tau_2=\alpha$ where $1+(\,(\alpha+1)^{1/2}+2)^{3/4}=0:$ \\

\noindent{\sf Step 2.1:} We  have that $\bar{n}_{11}=1$ and $\bar{n}_{12}=1$  and $\bar{n}_{13}=0$. \\

\noindent{\sf Step 2.2:} We have that
$$\begin{array}{ll}
   a_{1}=\limit_{t\rightarrow 2}\,\frac{p_2(t)}{p_{1}(t)}=(\,(\alpha+1)^{1/2}+2)^{1/4} \alpha &\\
  a_{0}=\limit_{t\rightarrow 2}\, p_1(t)f_1(t)=0,\quad &f_1(t):=\frac{p_2(t)}{p_{1}(t)^{2}}-a_{1}. \\
 \end{array}$$

\noindent{\sf Step 2.3:} We  obtain the asymptote $\widetilde{\cal C}_2$, defined by the proper parametrization
$$\widetilde{\cal Q}_2(t)=\left(t,\, ((\,(\alpha+1)^{1/2}+2)^{1/4} \alpha)t,\, p_3(\alpha)\right)=\left(t,\, ((\,(\alpha+1)^{1/2}+2)^{1/4} \alpha)t,\frac{\alpha+3}{\sin(\alpha)}\right),$$
 where $1+(\,(\alpha+1)^{1/2}+2)^{3/4}=0$.

\end{itemize}

\para

\noindent{\sf Step 3:} We observe that all the roots of $p_{12}(s)$ are also roots of $p_{22}(s)$ and $p_{23}(s)$ (that is, the poles of $p_1$ are the poles of $p_2$ and $p_3$), so there are no more  asymptotes.

\para

\noindent{\sf Step 4:} The algorithm {\sf returns} the   asymptotes  $\widetilde{\cal C}_i,\,i=1,2$ of the input curve, $\cal C$


\end{example}

\section{Conclusion}\label{S-conclusion}

The main results of this paper, Theorems \ref{T-final-ramas-G} and \ref{T-final-G}, provides a way to determine the branches and the generalized asymptotes of a curve parametrized not necessarily by two rational functions by only computing some simple limits of functions constructed from the given parametrization  defined by two   meromorphic functions on the complex plane. We remind that a meromorphic function on the complex plane is a function that is holomorphic on all of $\Bbb C$ except for a set of isolated points, which are poles of the function. We prove this theorem and we develop an efficient algorithm which determines all the branches and all the g-asymptotes which are obtained from the poles of the functions. As a complement, some corollaries are derived that allow us to obtain the horizontal and vertical asymptotes in an extremely simple way. This technique is proved to work on several illustrative examples.

\para

It is important to stress that this procedure can be trivially applied for dealing with   parametrizations of curves in  $n$--dimensional space. Thus, the present paper yields a remarkable improvement of the methodology developed in \cite{paper3} and \cite{newBP}.

\para

As a future work, we aim to extend the notion of g-asymptote to the study of the asymptotic behavior of algebraic surfaces. We look for surfaces which approach a given one of higher degree, when ``moving to infinity'', that is, when some of the coordinates take infinitely large values. The ideas introduced in this paper might provide the foundations for efficient methods that allow us to compute those ``asymptotic surfaces''.

\section*{Acknowledgements}
The author S. P\'erez-D\'{\i}az is partially supported by  Ministerio de Ciencia, Innovaci\'on y Universidades - Agencia Estatal de Investigaci\'on/PID2020-113192GB-I00 (Mathematical Visualization: Foundations, Algorithms and Applications).  The author R. Magdalena Benedicto is partially supported by the State Plan for Scientific and Technical Research and Innovation of the Spanish MCI (PID2021-127946OB-I00). \\
 The author S.\,P\'erez-D\'{\i}az belongs to the Research Group ASYNACS (Ref.CCEE2011/R34).\\

\noindent
{\bf Conflict of interest:} The authors declare no conflict of interest.\\

\noindent
{\bf Author contributions:} Authors contributed equally to this work and they worked together through the whole paper. All authors have read and agreed to the published version of the~manuscript.


\end{document}